\documentclass[12pt,oneside,reqno]{amsart}
\usepackage{amsmath,amssymb,amsthm,textcomp}
\usepackage{amsfonts,graphicx}
\usepackage[mathscr]{eucal}
\usepackage{color}
\usepackage[table]{xcolor}
\usepackage{booktabs}

\usepackage{subfig}
\usepackage{epsfig}
\usepackage{caption}
\theoremstyle{definition}
\usepackage{float}
\usepackage{mathtools} 
\usepackage[colorlinks=true,urlcolor=blue,citecolor=blue,linkcolor=blue,bookmarks=true]{hyperref} 
\usepackage[sort & compress,comma,authoryear]{natbib}
\floatstyle{plaintop}
\numberwithin{equation}{section}

\newcommand{\ncom}{\newcommand}

\ncom{\beq}{\begin{equation}}
	\ncom{\eeq}{\end{equation}}
\ncom{\bea}{\begin{eqnarray*}}
	\ncom{\eea}{\end{eqnarray*}}
\ncom{\beqa}{\begin{eqnarray}}
	\ncom{\eeqa}{\end{eqnarray}}
\ncom{\nno}{\nonumber}
\ncom{\non}{\nonumber}
\ncom{\ds}{\displaystyle}
\ncom{\half}{\frac{1}{2}}
\ncom{\mbx}{\makebox{.25cm}}
\ncom{\hs}{\mbox{\hspace{.25cm}}}
\ncom{\rar}{\rightarrow}
\ncom{\Rar}{\Rightarrow}
\ncom{\noin}{\noindent}
\ncom{\bc}{\begin{center}}
	\ncom{\ec}{\end{center}}
\ncom{\sz}{\scriptsize}
\ncom{\rf}{\ref}
\ncom{\s}{\sqrt{2}}
\ncom{\sgm}{\sigma}
\ncom{\Sgm}{\Sigma}
\ncom{\psgm}{\sigma^{\prime}}
\ncom{\dt}{\delta}
\ncom{\Dt}{\Delta}
\ncom{\lmd}{\lambda}
\ncom{\Lmd}{\Lambda}
\ncom{\Th}{\Theta}
\ncom{\e}{\eta}
\ncom{\eps}{\epsilon}
\ncom{\pcc}{\stackrel{P}{>}}
\ncom{\lp}{\stackrel{L_{p}}{>}}
\ncom{\dist}{{\rm\,dist}}
\ncom{\sspan}{{\rm\,span}}
\ncom{\re}{{\rm Re\,}}
\ncom{\im}{{\rm Im\,}}
\ncom{\sgn}{{\rm sgn\,}}
\ncom{\ba}{\begin{array}}
	\ncom{\ea}{\end{array}}
\ncom{\hone}{\mbox{\hspace{1em}}}
\ncom{\htwo}{\mbox{\hspace{2em}}}
\ncom{\hthree}{\mbox{\hspace{3em}}}
\ncom{\hfour}{\mbox{\hspace{4em}}}
\ncom{\vone}{\vskip 2ex}
\ncom{\vtwo}{\vskip 4ex}
\ncom{\vonee}{\vskip 1.5ex}
\ncom{\vthree}{\vskip 6ex}
\ncom{\vfour}{\vspace*{8ex}}
\ncom{\norm}{\|\;\;\|}
\ncom{\integ}[4]{\int_{#1}^{#2}\,{#3}\,d{#4}}
\ncom{\vspan}[1]{{{\rm\,span}\{ #1 \}}}
\ncom{\dm}[1]{ {\displaystyle{#1} } }
\ncom{\ri}[1]{{#1} \index{#1}}

\newtheorem{theorem}{\bf Theorem}[section]

\newtheorem{proposition}{Proposition}[section]
\newtheorem{lemma}{Lemma}[section]

\newtheoremstyle
{remarkstyle}
{}
{11pt}
{}
{}
{\bfseries}
{:}
{     }
{\thmname{#1} \thmnumber{#2} }

\theoremstyle{remarkstyle}

\def\eps{\varepsilon}

\begin{document}
	\title{\large R\lowercase{ecord-based}  \lowercase{transmuted}  \lowercase{unit~omega} \lowercase{distribution}: \lowercase{different~ methods~ of estimation} \lowercase{and} \lowercase{applications}} 
	\author[Ashok Kumar Pathak]{Ashok Kumar Pathak$^{1}$}
	\author{Mohd. Arshad$^{2}$}
	\author{Alok Kumar Pandey$^{1}$}
	\author{Alam Ali$^{1*}$
		\\
		$^{1}$D\lowercase{epartment of} M\lowercase{athematics and} S\lowercase{tatistics}, C\lowercase{entral} U\lowercase{niversity of} P\lowercase{unjab},\\
		B\lowercase{athinda}, P\lowercase{unjab}-151401, I\lowercase{ndia}.
		\\$^{2}$D\lowercase{epartment of} M\lowercase{athematics}, I\lowercase{ndian} I\lowercase{nstitute of} T\lowercase{echnology}
		I\lowercase{ndore},\\ S\lowercase{imrol}, I\lowercase{ndore}-453552, I\lowercase{ndia}.}
	\thanks{*Corresponding Author E-mail Address: alam2ali1996@gmail.com (Alam Ali)}
	\thanks{The research of  Alam Ali was supported by DST, Government of India.}
	
	\begin{abstract} 
		\cite{dombi2019omega} introduced a three parameter omega distribution and showed that its asymptotic distribution is the Weibull model. We propose a new record-based transmuted generalization of the unit omega distribution by considering \cite{balakrishnan2021record} approach. We call it the RTUOMG distribution. We derive expressions for some statistical quantities, like, probability density function, distribution, hazard function, quantile function, moments, incomplete moments, inverted moments, moment generating function, Lorenz curve, and Bonferroni curve of the proposed distribution. The numerical values of various measures of central tendency and coefficient of skewness and  kurtosis are also presented. Concepts of stochastic ordering and some results related to ordered statistics of the RTUOMG distribution are discussed. The parameters of the RTUOMG distribution are estimated using five distinct estimators. Additionally, the Monte Carlo simulations are performed to assess the performance of these estimators. Finally, two real data sets are analyzed to demonstrate the utility of the RTUOMG distribution.

	\end{abstract}
	
	\maketitle

%
%
	\section{Introduction}
	\noindent Proportion data is frequently encounter in diverse disciplines such as economics, finance, reliability, medical, biology, chemistry, etc. It is known that when  the values of a random variable are reported in percentage or fraction of the whole, it is referred to as proportion data. Such variables encompass all possible values within the unit interval and require an appropriate probability distribution to effectively model the observed data. Modeling approaches on the bounded interval have received considerable attention  since they are related to specific issues such as the acceptance rate, recovery rate, mortality rate, scores, proportion of the educational measurements, etc. In the last decade, a lot of research have  been done for developing distributions that are defined on a bounded interval. Various techniques are used to generate these distributions such as random variable transformation, function composition, and generation of new families of distributions.
	
	Beta distribution is one of the popularly known examples of the bounded distribution which has different density shapes and is widely used in different areas of applied sciences. Several alternative models have been defined on the bounded interval in the statistical distribution literature in order to provide better results. Topp-Leone and Kumaraswamy distributions are also well known example of the unit distribution introduced by \cite{topp1955family} and \cite{kumaraswamy1980generalized}, respectively. These models have received great attention in the past years in diverse areas of the applied sciences. However,  these models do not have a closed expression for the moments and are mathematically less tractable. Due to growing interest in study of the bounded distributions, some new families of the unit distributions have been proposed and studied in the recent literature. Some new among these are, namely,  unit-inverse gaussian distribution \cite{ghitany2019unit}, unit-Lindley distribution \cite{mazucheli2019one},  unit Weibull \cite{mazucheli2020unit},  unit log-log distribution \cite{ribeiro2021unit}, unit-bimodal Birnbaum-Saunders \cite{martinez2022unit}, unit Muth distribution \cite{maya2022unit}, \cite{korkmaz2023unit}, and  unit generalized half-normal distribution \cite{mazucheli2023unit}.
	
	Recently, by means of the omega function, \cite{dombi2019omega} proposed a new probability distribution on the bounded domain and discussed its applications in reliability theory. A random variable $U$ is said to have omega distribution with parameters $\alpha$, $\beta$, and $d$, denoted by $\text{OMG}(\alpha, \beta, d)$ if its density function is given by
	\begin{equation}
		f_{U}(x; \alpha,\beta, d)=2\alpha\beta x^{\beta-1}\frac{d^{2\beta}}{d^{2\beta}-x^{2\beta}} \omega_{d}^{(-2\alpha,\beta)}(x),\; 0<x<d,
	\end{equation}
	where $\alpha>0$, $\beta>0$, $d>0$, and $\omega_{d}^{(\alpha,\beta)}(x)= \left(\frac{d^\beta+x^\beta}{d^\beta-x^\beta}\right)^{\alpha d^\beta/2}$ denotes the omega function.
	The hazard function of the omega distribution can be monotonic, constant, and bathtub shapes. Therefore, it can used in modeling diverse class of real phenomena.  \cite{dombi2019omega} showed that the limiting omega distribution is just the Weibull distribution and can be used in place of the Weibull distribution. For some more developments and discussion about omega distribution one can refer to \cite{okorie2019omega}, \cite{alsubie2021omega}, \cite{ozbilen2022bivariate}, \cite{birbiccer2022parameter}, and \cite{jonas2022generalized}. Specially, when $d=1$, the omega distribution corresponds to the unit omega distribution supported on the unit interval. Its density and distribution functions (DF) are given by
	\begin{equation}
		f_{U}(x; \alpha,\beta)=2\alpha\beta x^{\beta-1}\frac{1}{1-x^{2\beta}} \omega_{1}^{(-2\alpha,\beta)}(x),\; 0<x<1
	\end{equation}
	and 
	\begin{equation}
		F_{U}(x; \alpha,\beta)=1- \omega_{1}^{(-2\alpha,\beta)}(x),\; 0<x<1,
	\end{equation}
	respectively. \cite{birbiccer2022parameter} showed that the unit omega distribution is actually a unit exponentiated half logistic distribution. That is, the random variable $X=e^{-Z}$ has the unit omega distribution  when $Z$ is the exponentiated half logistic random variable. The distribution function of the unit omega is mathematically simple  and independent special functions and  anticipate ease in exploring its important statistical properties and statistical inference in  the model. It can have different shapes including U-shaped, J-shaped, reversed J-shaped, left and right skewed and are used in various applications related to reliability.
	
	By the means of the quadratic rank transformation, \cite{shaw2009alchemy} introduced a new class of transmuted distributions using baseline distribution. We say that a random variable $Y$ has the transmuted distribution if its distribution is expressed as
	\begin{equation}\label{TD1}
		{G}_{\theta}(y)=(1+\theta)H(y)-\theta H^2(y),\;y\in \mathbb{R},
	\end{equation} 
	where $H(\cdot)$ is the baseline distribution and $|\theta|\leq 1$. One can see that the transmuted distribution defined in (\ref{TD1}) can be obtained from the combination of distribution functions of the smallest and largest order statistics from the baseline distribution in a sample of size 2. Transmutation of the baseline distribution is a powerful tool to construct the skewed probability distribution and have been used by several researchers in the recent years. \cite{granzotto2017cubic} proposed the cubic rank transmuted distribution and discussed its important statistical properties. In addition to these publications, there are several papers in the literature that discussed the fundamental characteristics of various new quadratic and cubic transmuted distributions, readers can see
	\cite{elbatal2013transmuted},  \cite{merovci2013transmuted}, \cite{tian2014transmuted},
	\cite{khan2014new}, \cite{kemaloglu2017transmuted}, \cite{alizadeh2018transmuted},
	\cite{kharazmi2021informational}, \cite{chhetri2022cubic}, and
	\cite{tanics2023cubic}. 
	Recently, \cite{balakrishnan2021record}  formulated a record-based transmuted map to generate
	new class of probability models from the baseline distribution. They showed that the record-based transmuted distribution are expressed via the relation
	\begin{equation}
		F_{p}(x)=H(x)+p\overline{H}(x)\log \overline{H}(x), \;x\in \mathbb{R},
	\end{equation}
	where $0\leq p\leq 1$, and $H(\cdot)$, $\overline{H}(\cdot)$ denotes the distribution function (DF) and survival function (SF) of the baseline distribution, respectively. \cite{balakrishnan2021record} also introduced some new record-based transmuted (RT) probability distribution, namely, RT-exponential distribution, RT-Weibull distribution, and RT-Linear exponential distribution. \cite{tanics2022record} discussed the different estimation  for parameters estimation of the RT-Weibull distribution and demonstrated its real application. The  RT-generalized linear exponential distribution is proposed and studied by \cite{arshad2022record}. They also analyzed the lifetime data sets using it.
	
	This paper aims to introduce a new record-based transmuted generalization of the unit omega distribution and study its important properties. The proposed distribution is a more broad family of probability distributions and includes the unit omega distribution as a submodel. Its density and hazard function can be used to handle a vast class of data that arises in a variety of domains. These functions have various shapes, including as increasing, decreasing, and bathtub shapes. 
	
	The organization of the article is as follows: Section 2 presents the mathematical framework of the proposed RTUOMG distribution. We derive the expressions for the DF, density, hazard function, quantile function, moments,  incomplete moments, inverted moments, moment generating function, lorenz curve, and bonferroni curve of the RTUOMG distribution. Section 3 presents some concepts of stochastic ordering and some results related to ordered statistic of the proposed distribution.  In Section 4, the maximum likelihood estimators, least squares and weighted least squares estimators, Cram{\'e}r-von Mises estimators, and Anderson-Darling estimators of the parameters are explored.  Section
	5 presents the detailed Monte Carlo simulation study to validate the performances of the
	estimators through the absolute bias and mean squared error measures. Finally, in Section 6, two real data sets are analyzed to show the applicability of RTUOMG distribution in real situations. 
	\section{Record-based transmuted unit omega distribution} \noindent Let $Y_1,Y_2,\ldots$ be a sequence of independent and identically distributed (iid) random variables having DF $H(\cdot)$. Let $Y_{U(1)}$ and $Y_{U(2)}$ denote the first two upper records from the sequence of iid random variables. Now, define a random variable $X$ as
	\begin{equation*}
		X=\begin{cases}
			Y_{U(1)},\;\text{with probability}~ 1-p\\
			Y_{U(2)},\;\text{with probability}~ p,
		\end{cases}
	\end{equation*}
	where $0\leq p\leq 1$. The DF of $X$ is obtained as follows (see \cite{balakrishnan2021record}):
	\begin{align}
		F_{X}(x)=& (1-p) P(Y_{U(1)}\leq x)+ p P(Y_{U(2)}
		\leq x)\nonumber\\
		=& (1-p) H(x)+p\left[1-\overline{H}(x)\sum_{r=0}^{1}\frac{(-\log  \overline{H}(x)
			)^{r}}{r!}\right]\nonumber\\
		=& (1-p)H(x)+p[1- \overline{H}(x)(1-\log\overline{H}(x))]\nonumber\\
		=& (1-p)H(x)+p[H(x)+\overline{H}(x) \log \overline{H}(x)]\nonumber\\
		=& H(x)+p \overline{H}(x) \log \overline{H}(x), \;x\in \mathbb{R}, 
	\end{align}
	where $ \overline{H}(x)=1-H(x)$ is SF of the baseline distribution $H(x)$. The density function and hazard function (HF) of $X$ are respectively, given by 
	\begin{equation}
		f_{X}(x)=h(x)[1-p-p\log \overline{H}(x)], \;x\in \mathbb{R}
	\end{equation}
	and 
	\begin{equation}
		r_{X}(x)=r(x)\frac{1-p-p\log \overline{H}(x)}{1-p\log \overline{H}(x)}, \;\;x\in \mathbb{R},
	\end{equation}
	where $h(x)$ is the density function of the baseline distribution, and $r(x)$ denotes the HF of the baseline distribution and is defined as $r(x)=h(x)/ \overline{H}(x)$.
	
\begin{figure}[h]
		\includegraphics[scale=0.50, angle=0]{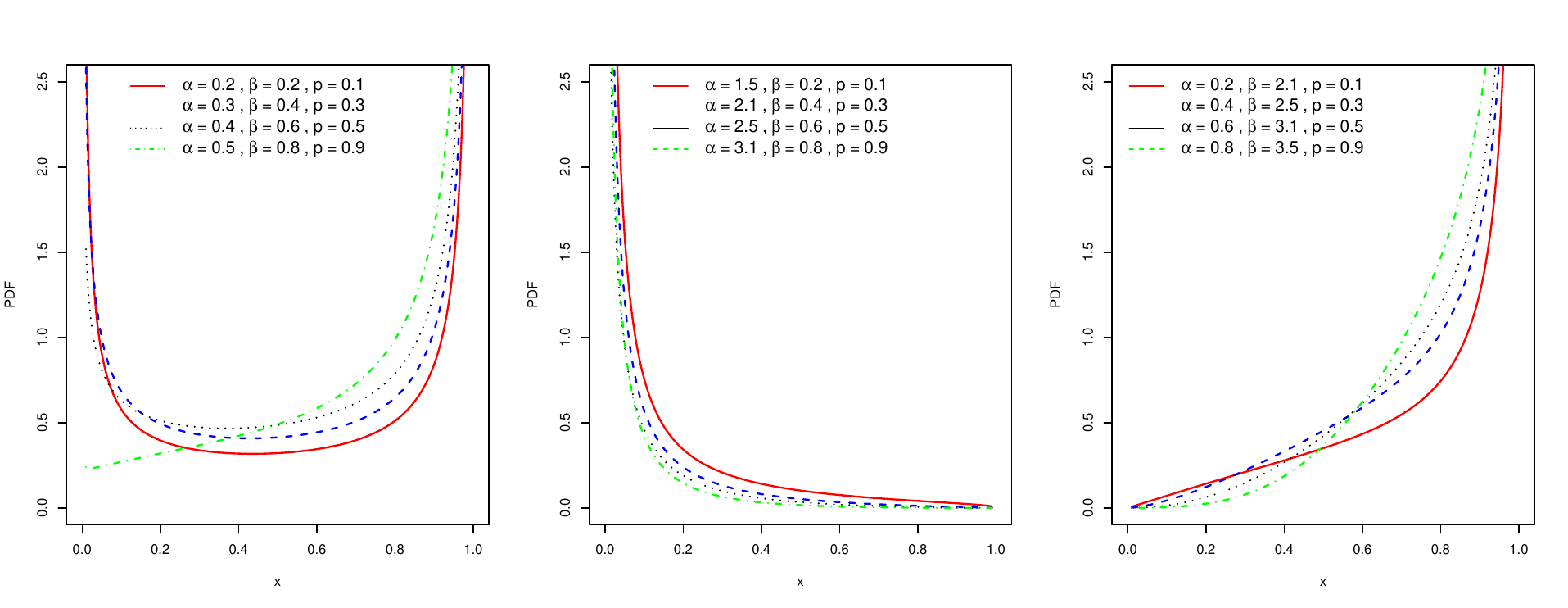}
		\caption{PDF plots of RTUOMG distribution.}
		\label{density}
	\end{figure}

\begin{figure}[h]
		\includegraphics[scale=0.450, angle=0]{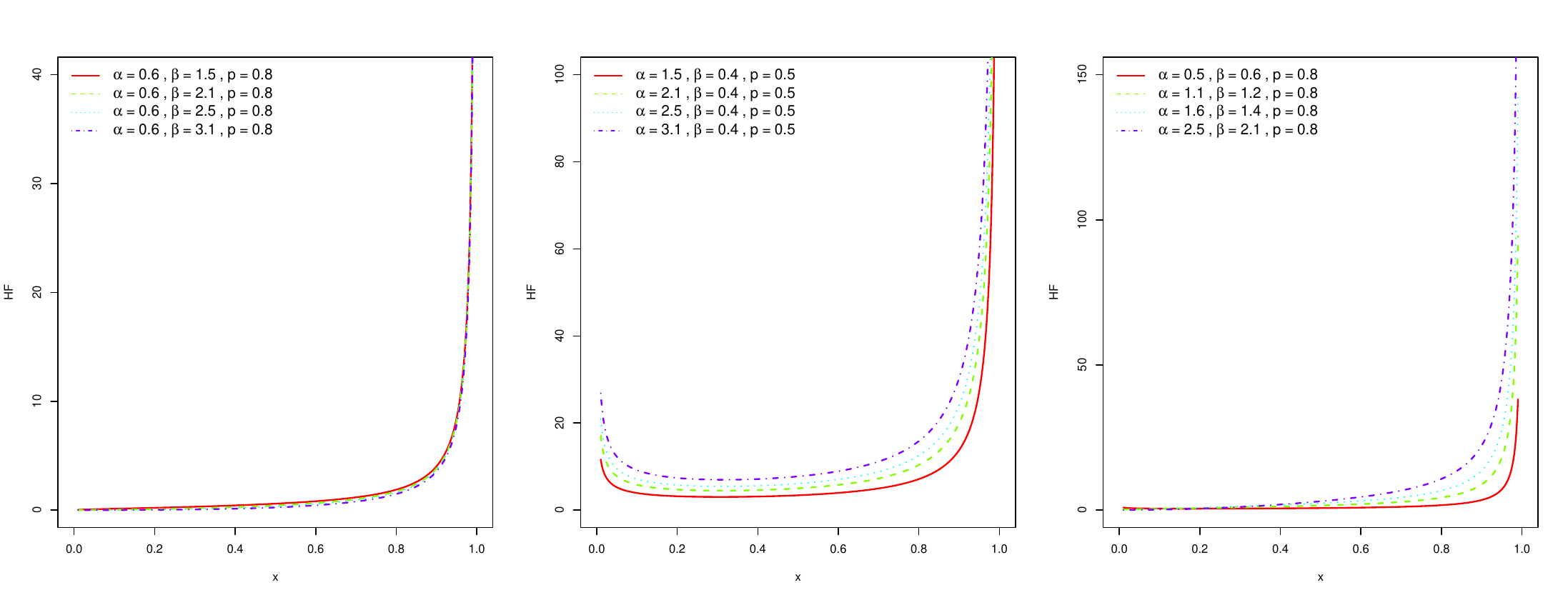}
		\caption{HF plots of RTUOMG distribution.}
		\label{fighazard}
	\end{figure}

	We say that a random variable $X$ follows the record-based transmuted unit omega distribution (denoted by $X\sim \text{RTUOMG}(\alpha, \beta, p)$) if its distribution function is given by 
	\begin{equation}\label{CDF}
		F_{X}(x;\alpha,\beta,p)=1- \omega_{1}^{(-2\alpha,\beta)}(x)+p\omega_{1}^{(-2\alpha,\beta)}(x)\log \left(\omega_{1}^{(-2\alpha,\beta)}(x)\right),  
	\end{equation}
	where $0<x<1$, $\alpha>0$, $\beta>0$ and $p\in [0,1]$. The corresponding density function of the RTUOMG distribution is given by 
	\begin{equation}\label{pdf1}
		f_{X}(x;\alpha,\beta,p)=\left(\frac{2\alpha\beta x^{\beta-1}}{1-x^{2\beta}}\right)\omega_{1}^{(-2\alpha,\beta)}(x)\left[ 1-p-p \log \left(\omega_{1}^{(-2\alpha,\beta)}(x)\right) \right].
	\end{equation}
	It may be seen that when $p=0$, the RTUOMG distribution reduces to OMG distribution. Also, one can easily show that $F_{X}(x)\leq H(x), \; \forall x\in (0,1)$. The HF of the the RTUOMG distribution is 
	\begin{equation}
		r_{X}(x;\alpha,\beta,p)=\frac{2\alpha\beta x^{\beta-1}}{(1-x^{2\beta})}\left\{\frac{1-p-p \log \left(\omega_{1}^{(-2\alpha,\beta)}(x)\right)}{1+p \log \left(\omega_{1}^{(-2\alpha,\beta)}(x)\right)}\right\}.
	\end{equation}
	For different values of the model parameters $\alpha$, $\beta$, and $p$, the plots of the probability density function (PDF) and HF are shown in the Figure \ref{density} and Figure \ref{fighazard}, respectively. From these figures, we see that the PDF and HF of the  RTUOMG distribution takes different shapes and suggest applicability of the model in diverse areas.
	\subsection{Quantile function}
	The quantile function of the RTUOMG distribution can be derived in terms of the Lambert W function which is defined as 
	\begin{equation*}
		W(\xi) e^{W(\xi)}=\xi,
	\end{equation*}
	where $\xi$ is a complex number. When $\xi \geq-1/e$, it has two branches, namely, principal branch ($W_{o}$)  and negative branch ($W_{-1}$). For a detailed discussion about the Lambert W function one can refer to a recent paper by \cite{jodra2022intermediate}.
	\begin{theorem}
		Let X be a random variable having RTUOMG$(\alpha,\beta,p)$ distribution. Then, the quantile function $Q(u)$ is
		\begin{equation}\label{eq:qf}
			Q(u)= \left[\frac{e^{\frac{-1}{\alpha}\left\{ W_{-1}\left(\frac{u-1}{pe^{1/p}}\right)+\frac{1}{p}\right\}}-1}{e^{\frac{-1}{\alpha}\left\{ W_{-1}\left(\frac{u-1}{pe^{1/p}}\right)+\frac{1}{p}\right\}}+1}\right]^{\frac{1}{\beta}},~ 0<u<1.
		\end{equation}
	\end{theorem}
	\begin{proof}
		For $u\in (0,1)$, the solution of  $F(x)=u$ yields the quantile function. Now, we have 
		\begin{align*}
			\left\{1-p\log\left(\frac{1+x^\beta}{1-x^\beta}\right)^{-\alpha}\right\} \left(\frac{1+x^\beta}{1-x^\beta}\right)^{-\alpha}=& 1-u\\
			\left\{ \frac{-1}{p} +\log\left(\frac{1+x^\beta}{1-x^\beta}\right)^{-\alpha}\right\} e^{\log\left(\frac{1+x^\beta}{1-x^\beta}\right)^{-\alpha}-\frac{1}{p}}=& \frac{u-1}{p e^{\frac{1}{p}}}.
		\end{align*}
		It can be verified that $ (u-1)/pe^{1/p}\in [-1/e,0)$, and $\frac{-1}{p} +\log\left(\frac{1+x^\beta}{1-x^\beta}\right)^{-\alpha} \in (-\infty,-1].$ Now, using the $W_{-1}$ function in the above equation, we get
		\begin{equation}\label{eq:1111}
			W_{-1}\left(\frac{u-1}{pe^{\frac{1}{p}}}\right)= \frac{-1}{p} +\log\left(\frac{1+x^\beta}{1-x^\beta}\right)^{-\alpha}.
		\end{equation}
		On solving \eqref{eq:1111} for $x$, completes the proof of the theorem.
	\end{proof}
	Next, we have the following lemmas which are needed to derive various types of the moments.
	\begin{lemma}\label{lemma1} For $\lambda>0$ and $\mu>0$, we have (see [\cite{lee2011generalization}])
		\begin{equation} \label{eq:lemma}
			\int_{0}^{1} x^{\lambda-1} (1-x)^{\mu-1} (1-\beta x)^{-\nu} dx= B(\lambda,\mu)~ _{2}F_{1}(\nu,\lambda;\lambda+\mu;\beta),
		\end{equation}
		where $B(\alpha,\beta)$ and $_{2}F_{1}(\alpha,\beta;\gamma;x)$ are the beta function and Gauss hypergeometric function, respectively, defined by
		\begin{equation*}
			B(\alpha,\beta)= \int_{0}^{1} t^{\alpha-1} (1-t)^{\beta-1} dt
		\end{equation*}
		and
		\begin{equation*}
			_{2}F_{1}(\alpha,\beta;\gamma;x)= \sum_{i=0}^{\infty} \frac{(\alpha)_{i} (\beta)_{i}}{(\gamma)_{i}} \frac{x^i}{i!},
		\end{equation*}
		where $(n)_{i}=n(n+1)(n+2)\ldots(n+i-1)$ denotes the falling factorial.
	\end{lemma}
	\begin{lemma}[Formula 1.513.1 in \cite{gradshteyn2014table}]\label{lemma2}
		If $y^2<1$, then
		\begin{equation}\label{eq:ser}
			\log\left(\frac{1+y}{1-y}\right)= 2\sum_{k=1}^{\infty} \frac{1}{2k-1} y^{2k-1}.
		\end{equation}
	\end{lemma}
	\begin{lemma} \label{lemma3}
		For $Re(a)>Re(b)>0, |arg(1-x)|<\pi$, we have (see \cite{ozarslan2019some})
		\begin{equation} 
			\int_{0}^{\lambda} y^{b-1} (1-y)^{a-b-1} (1-xy)^{-\alpha} dy= B(b,a-b)~ _{2}\delta_{1}(\alpha,[b,a;\lambda],x),
		\end{equation}
		where $B(\alpha,\beta)$ is the beta function defined in Lemma \ref{lemma1} and $_{2}\delta_{1}(\alpha,[b,a;\lambda],x)$ represents the incomplete gauss hypergeometric function defined by
		\begin{equation*}
			_{2}\delta_{1}(\alpha,[b,a;\lambda],x)= \sum_{i=0}^{\infty} (\alpha)_{i} [b,a;\lambda]_{i} \frac{x^i}{i!},
		\end{equation*}
		where $[b,a;\lambda]_{i}$ denotes the incomplete Pochhammer ratio which is introduced in terms of the incomplete beta function as follows
		\begin{equation*}
			[b,a;\lambda]_{i}= \frac{B_{\lambda}(b+i,a-b)}{B(b,a-b)},
		\end{equation*}
		where $B_{\lambda}(m,n)$ is known as incomplete beta function and given by
		\begin{equation*}
			B_{\lambda}(m,n)= \int_{0}^{\lambda} t^{m-1} (1-t)^{n-1} dt, ~ Re(m)>0,~Re(n)>0,~ 0\leq \lambda<1.
		\end{equation*}
	\end{lemma}
	\subsection{Moments and measures} The $r$th order moment is defined by
	\begin{equation*}
		\mu'_{r}=\mathbb{E}[X^r]=\int x^r f_{X}(x)dx.
	\end{equation*}
	\begin{proposition}\label{proposition1}
		Let $X\sim \text{RTUOMG}(\alpha, \beta, p)$. Then
		\begin{equation}\label{eq:moment}
			\begin{aligned}
				\mu'_{r}=& 2\alpha(1-p) B\Big(\frac{r}{\beta}+1, \alpha\Big)~{}_{2}F_{1} \Big(\alpha+1,\frac{r}{\beta}+1; \frac{r}{\beta}+\alpha+1;-1\Big) \\
				& + 4p\alpha^2 \sum_{k=1}^{\infty}\left[\frac{1}{2k-1} B\Big(\frac{r}{\beta}+2k, \alpha\Big)~{}_{2}F_{1} \Big(\alpha+1,\frac{r}{\beta}+2k; \frac{r}{\beta}+2k+\alpha;-1\Big)\right].		
			\end{aligned}
		\end{equation}
	\end{proposition}
	\begin{proof}
		We have
		\begin{align*}
			\mu'_{r}=&\int_{0}^{1}x^{r}	\left(\frac{2\alpha\beta x^{\beta-1}}{1-x^{2\beta}}\right)\omega_{1}^{(-2\alpha,\beta)}(x)\left[ 1-p-p \log \left(\omega_{1}^{(-2\alpha,\beta)}(x)\right) \right]dx\\
			=&\int_{0}^{1} x^r \frac{2\alpha\beta x^{\beta-1}}{1-x^{2\beta}} \left(\frac{1+x^\beta}{1-x^\beta}\right)^{-\alpha} \left[ 1-p+p\alpha \log\left(\frac{1+x^\beta}{1-x^\beta}\right)  \right] dx.
		\end{align*}
		Setting $y=x^{\beta}$, we get
		\begin{equation*}\label{eq1}
			\mu'_{r}= \int_{0}^{1} y^{\frac{r}{\beta}} \frac{2\alpha}{1-y^2} \left(\frac{1+y}{1-y}\right)^{-\alpha} \left[1-p+p\alpha \log\left(\frac{1+y}{1-y}\right) \right] dy.	
		\end{equation*}
		Using Lemma (\ref{lemma2}), we get
		\begin{equation*}
			\mu'_{r}= 2\alpha \int_{0}^{1}  \frac{y^{\frac{r}{\beta}}}{(1+y)(1-y)} \left(\frac{1+y}{1-y}\right)^{-\alpha} \left[1-p+2p\alpha \sum_{k=1}^{\infty} \frac{1}{2k-1} y^{2k-1} \right] dy.
		\end{equation*}
		By an use of the Lemma (\ref{lemma1}) and simple calculation completes the proof.
	\end{proof}
	\noindent The coefficient of Skewness (CS) and the coefficient of Kurtosis (CK) using first four moment are respectively, defined as
	\begin{equation*}
		CS= \frac{\mu'_{3}-3\mu'_{2}\mu'_{1}+2{\mu'_{1}}^3}{{\mu_{2}}^{\frac{3}{2}}}
	\end{equation*}
	and
	\begin{equation*}
		CK= \frac{\mu'_{4}-4\mu'_{1}\mu'_{3}+6{\mu'_{1}}^2\mu'_{2}-3{\mu'_{1}}^4}{{\mu_{2}}^2},
	\end{equation*}
	where ${\mu_{2}}$ is variance and is defined as $	\mu_{2}= \mu'_{2}- (\mu'_{1})^{2}$.
	
\begin{figure}[t] 
		\centering
		\subfloat{\includegraphics[scale=0.53, angle=0]{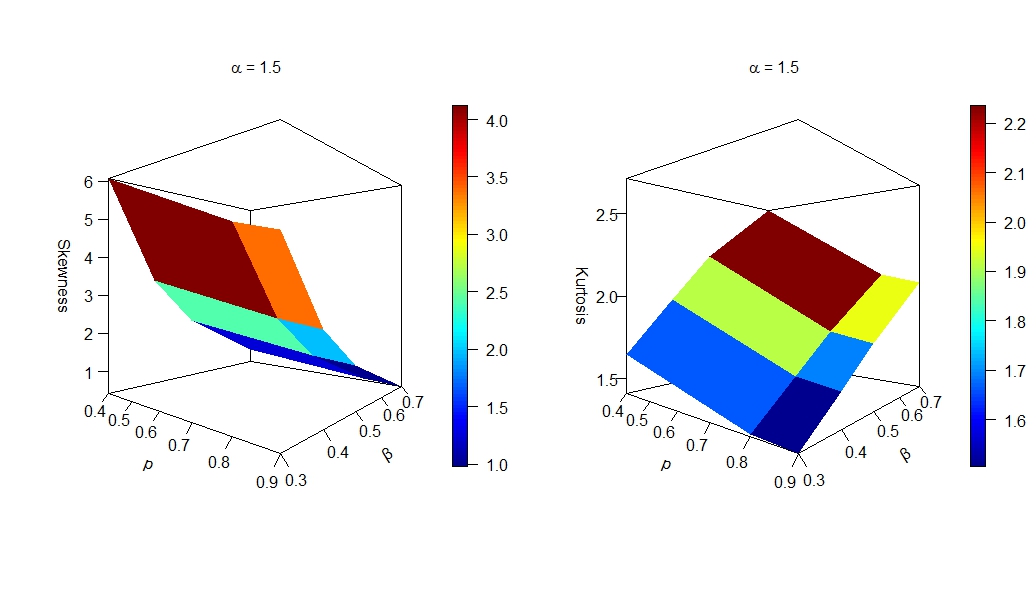}}\\
		\subfloat{\includegraphics[scale=0.53, angle=0]{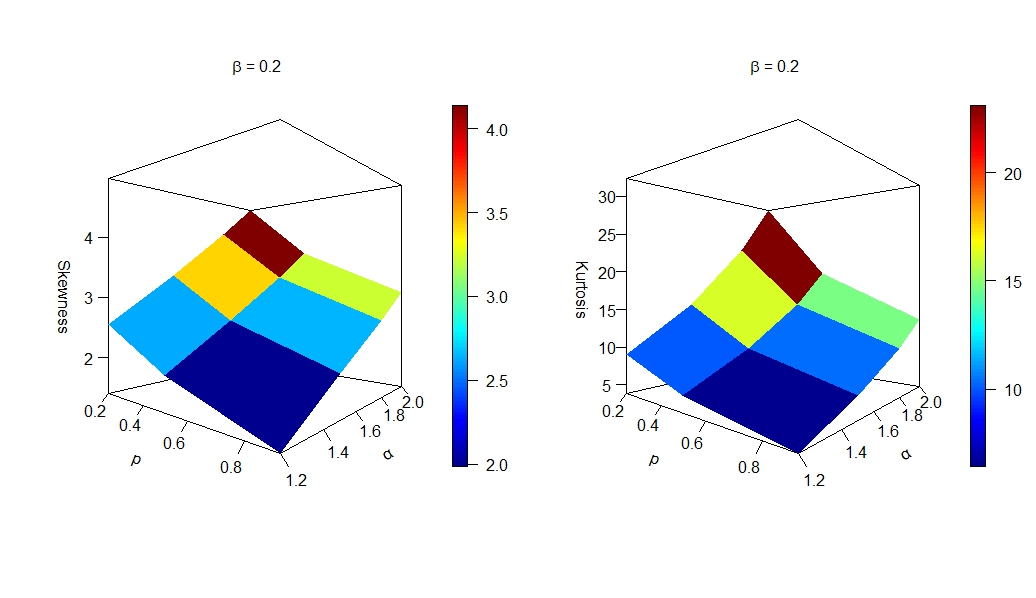}}\\
		\caption{ Graph of Skewness and Kurtosis for different parameters.}
		\label{fig.6}
	\end{figure}

	\noindent In Table \ref{T1}, we find a compilation of statistical moments, variance, CS, and CK for various parameter values of the RTUOMG$(\alpha,\beta,p)$ distribution.  Additionally, the Figure \ref{fig.6} provides a visual representation of how skewness and kurtosis change across different parameter configurations.
	
	\subsection{Moment generating function}
	Let $X\sim \text{RTUOMG}(\alpha,\beta,p)$. Then, the moment generating function  $M_{RTUOMG}(t)$ of $X$ is 
	\begin{equation}
		\begin{aligned}
			M_{RTUOMG}(t)&= \int_{0}^{1} e^{tx} f_X(x) dx= \sum_{r=0}^{\infty} \frac{t^r}{r!} \int_{0}^{1} x^r f_X(x) dx\\
			&=2\alpha(1-p)\sum_{r=0}^{\infty}\left\{\frac{t^r}{r!} B\Big(\frac{r}{\beta}+1, \alpha\Big)~{}_{2}F_{1} \Big(\alpha+1,\frac{r}{\beta}+1; \frac{r}{\beta}+1+\alpha;-1\Big)\right\} \\
			& + 4p\alpha^2 \sum_{r=0}^{\infty}\left\{\frac{t^r}{r!}\sum_{k=1}^{\infty}\left[\frac{1}{2k-1} B\Big(\frac{r}{\beta}+2k, \alpha\Big)~{}_{2}F_{1} \Big(\alpha+1,\frac{r}{\beta}+2k; \frac{r}{\beta}+2k+\alpha;-1\Big)\right]\right\}.
		\end{aligned}
	\end{equation}                                                                   Next, we present the incomplete moments for the RTUOMG distribution.                                                                                                                                                                                                                                                                                                                                                                                                                                                                                                                                                                                \subsection{Incomplete moment} The $r$th incomplete moment is defined as
	\begin{equation*}
		\phi_{r}(z)=\int_{0}^{z} x^r f_{X}(x)dx.
	\end{equation*}
	\begin{proposition}\label{proposition2}
		Let $X\sim\text{RTUOMG}(\alpha,\beta,p)$. Then
		\begin{equation}\label{incompletemoment}
			\begin{aligned}
				\phi_{r}(z)&= 2\alpha(1-p) B\Big(\frac{r}{\beta}+1, \alpha\Big)~{}_{2}\delta_{1} \Big(\alpha+1,\left[\frac{r}{\beta}+1, \alpha+\frac{r}{\beta}+1;z^\beta\right],-1\Big) \\
				& + 4p\alpha^2 \sum_{k=1}^{\infty}\left[\frac{1}{2k-1} B\Big(\frac{r}{\beta}+2k, \alpha\Big)~{}_{2}\delta_{1} \Big(\alpha+1,\left[\frac{r}{\beta}+2k, \alpha+\frac{r}{\beta}+2k;z^\beta\right],-1\Big)\right].
			\end{aligned}
		\end{equation}	
	\end{proposition}
	\noindent One can easily establish the proof of the Proposition \ref{proposition2}  with simple calculation and using the Lemma \eqref{lemma2}-\eqref{lemma3}.	
	
	\subsection{Inverted moments} The rth inverted moment of RTUOMG distribution  defined by 
	\begin{equation*}
		\mu^* _{r}= \int x^{-r} f_{X}(x)dx
	\end{equation*}
	\begin{proposition}\label{proposition3}
		Let $X\sim\text{RTUOMG}(\alpha,\beta,p)$. Then, the $r$th inverted moment is
		\begin{align}
			\mu^* _{r}&= 2\alpha(1-p) B\Big(1-\frac{r}{\beta}, \alpha\Big)~{}_{2}F_{1} \Big(\alpha+1,1-\frac{r}{\beta}; 1-\frac{r}{\beta}+\alpha;-1\Big) \nonumber\\
			& + 4p\alpha^2 \sum_{k=1}^{\infty}\left[\frac{1}{2k-1} B\Big(2k-\frac{r}{\beta}, \alpha\Big)~{}_{2}F_{1} \Big(\alpha+1,2k-\frac{r}{\beta}; 2k-\frac{r}{\beta}+\alpha;-1\Big)\right].
		\end{align}
	\end{proposition}
	\subsection{Lorenz and Bonferroni curve} The Lorenz curve measure the uncertainty in data. It is defined in terms of the incomplete moment and mean as 
	\begin{equation} \label{Lorenz}
		L(x)= \frac{\phi_{1}(z)}{\mathbb{E}(X)}.
	\end{equation}
	For the RTUOMG distribution, the expression of the Lorenz curve $L$ can be obtained with the help of (\ref{eq:moment}) and (\ref{incompletemoment}).\\
	The Bonferroni curve is obtained using the Lorenz curve as 
	\begin{equation}
		B(x)= \frac{L(x)}{F_X(x)}.
	\end{equation}
	For the RTUOMG distribution, we can calculate $B(x)$ by values of the its  component.
	\begin{table}[h]
		\caption{Mean, variance, coefficient of skewness (CS), and coefficient of kurtosis (CK) for different values of parameters $\alpha$, $\beta$, and $p$.} 
		\centering 
		\footnotesize
		\begin{tabular}{cccccccccc} 
			\hline
			$\beta$ & $p$& $\alpha$ &$\mu'_{1}$&$\mu'_{2}$&$\mu'_{3}$&$\mu'_{4}$&$\mu_{2}$ &CS&CK  \\
			\hline
			0.7& 0.2& 0.50& 0.5542510& 0.4285096& 0.3644273& 0.3239216& 0.1213154& -0.1787491&  1.5419750\\
			&& 1.25& 0.2901651& 0.1567514& 0.1047541& 0.0777581& 0.0725556&  0.8782454&  2.6772931\\
			&& 2.00 &0.1805073& 0.0711336& 0.0378631& 0.0235777& 0.0385507&  1.4672175&  4.6838320\\
			&& 2.75& 0.1247350& 0.0373275& 0.0160948& 0.0084739& 0.0217687&  1.8706299&  6.7571180\\
			&&3.50& 0.0923828& 0.0217244& 0.0077187& 0.0034674& 0.0131898&  2.1617793&  8.6740660\\
			&& 4.25& 0.0718269& 0.0136545& 0.0040625& 0.0015718& 0.0084954&  2.3771614& 10.355870\\
			&& 5.00& 0.0578708& 0.0091041& 0.0023016& 0.0007740& 0.0057551 & 2.5393587& 11.790823\\
			&&5.75& 0.0479085& 0.0063591& 0.0013840& 0.0004081& 0.0040639&  2.6634355& 12.997680\\
			&&6.50& 0.0405136& 0.0046111& 0.0008741& 0.0002279& 0.0029698&  2.7597558& 14.006293\\
			\hline
			$\alpha$ & $p$& $\beta$ &&&&&&&  \\
			1.5& 0.5& 0.50& 0.2310710& 0.1156176& 0.0734855& 0.0524786& 0.0622238&  1.1605455& 3.368963\\
			&& 1.25& 0.4610186& 0.2805956& 0.1949608& 0.1463780& 0.0680574&  0.1604957& 1.978038\\
			&& 2.00& 0.5881271& 0.4001152& 0.2959190& 0.2310710& 0.0542216& -0.2517045& 2.168528\\
			&& 2.75& 0.6672679& 0.4872244& 0.3768347& 0.3033692& 0.0419780& -0.4994430& 2.523453\\
			&& 3.50& 0.7210560& 0.5528266& 0.4420729& 0.3645541& 0.0329048& -0.6698945& 2.881894\\
			&& 4.25& 0.7599340& 0.6037974& 0.4954208& 0.4165299& 0.0262978& -0.7959642& 3.211055\\
			&& 5.00& 0.7893271& 0.6444595& 0.5397056& 0.4610186& 0.0214221& -0.8936171& 3.505502\\
			&& 5.75& 0.8123212& 0.6776175& 0.5769850& 0.4994245& 0.0177518& -0.9717678& 3.767074\\
			&& 6.50& 0.8307971& 0.7051561& 0.6087628& 0.5328590& 0.0149322& -1.0358663& 3.999452\\
			\hline 
			$\alpha$ & $\beta$& $p$ &&&&&&&\\
			0.5& 0.3& 0.1& 0.3710223& 0.2728817& 0.2271800& 0.1993273& 0.1352242 & 0.5146700& 1.679386\\
			&& 0.2& 0.4031998& 0.3038938& 0.2565056& 0.2271873& 0.1413237&  0.3766815& 1.533736\\
			&& 0.3& 0.4353719& 0.3349045& 0.2858282& 0.2550473& 0.1453558&  0.2427470& 1.437843\\
			&& 0.4& 0.4675487& 0.3659152& 0.3151508& 0.2829073& 0.1473134&  0.1116919& 1.386651\\
			&& 0.5& 0.4997255& 0.3969259& 0.3444733& 0.3107672& 0.1472004& -0.0177101& 1.377406\\
			&& 0.6& 0.5319022& 0.4279366& 0.3737959& 0.3386272& 0.1450166& -0.1465645& 1.409213\\
			&& 0.7& 0.5640815& 0.4589473& 0.4031185& 0.3664872& 0.1407593& -0.2757944& 1.482790\\
			&& 0.8 &0.5962576& 0.4899580& 0.4324410& 0.3943472& 0.1344348& -0.4060726& 1.600227\\
			&& 0.9& 0.6284332& 0.5209687& 0.4617636& 0.4222072& 0.1260404& -0.5373839& 1.764014\\
			\hline
		\end{tabular}\label{T1}
	\end{table}
	\newpage
	\section{Stochastic ordering and order statistics}
	\noindent In distribution theory and statistics, a stochastic ordering quantifies the concept of one random variable being bigger than another. A random variable $X$ is said to be smaller than a random variable $Y$ if
	\begin{enumerate}
		\item Stochastic order  $X\leq_{st}Y$ if $F_{X}(x)\geq F_{Y}(x)$ $\forall$ x,
		\item Hazard rate order $X\leq_{hr}Y$ if $h_{X}(x)\geq h_{Y}(x)$ $\forall$ x,
		\item Mean residual life order $X\leq_{mrl}Y$ if $m_{X}(x)\geq m_{Y}(x)$ $\forall$ x,
		\item Likelihood ratio order $X\leq_{lr}Y$ if $\frac{f_{X}(x)}{f_{Y}(x)}$ decreases in x.
	\end{enumerate} 
	\noindent  \cite{shaked2007stochastic} discussed the following interconnections\\
	$X\leq_{lr}Y \Rightarrow X\leq_{hr}Y \Rightarrow X\leq_{mrl}Y\Rightarrow X\leq_{st}Y $.\\
	Next, we have the following result.
	\begin{theorem}
		Let $X \sim$ RTUOMG$(\alpha_{1},\beta_{1},p_{1})$  and $Y \sim$ RTUOMG$(\alpha_{2},\beta_{2},p_{2})$. If $\alpha_{1}>\alpha_{2}$, $ \beta_{1}=\beta_{2}=\beta$, and  $p_{1}=p_{2}=p$,  then $X\leq_{lr}Y$. Hence, $X\leq_{hr}Y, X\leq_{mrl}Y$ and $X\leq_{st}Y$.
	\end{theorem}
	\textbf{Proof:} We have 
	\begin{equation*}\label{stochastic}
		\frac{f_{X}(x)}{f_{Y}(x)}= \frac{\frac{2\alpha_{1}\beta_{1} x^{\beta_{1}-1}}{1-x^{2\beta_{1}}} \left(\omega_{1}^{(-2\alpha_1,\beta_1)}(x)\right) \left\{ 1-p_{1}\left(  1+\log \left(\omega_{1}^{(-2\alpha_1,\beta_1)}(x)\right)\right)\right\}}{\frac{2\alpha_{2}\beta_{2} x^{\beta_{2}-1}}{1-x^{2\beta_{2}}} \left(\omega_{1}^{(-2\alpha_2,\beta_2)}(x)\right) \left\{1-p_{2}\left(   1+\log \left(\omega_{1}^{(-2\alpha_2,\beta_2)}(x)\right)\right)\right\}}.
	\end{equation*}
	Taking log and differentiating both sides with respect to $x$, we get 
	\begin{align*}
		\frac{d}{dx}\log\left\{\frac{f_{X}(x)}{f_{Y}(x)}\right\}&= \frac{\beta_{1}-1}{x} -\frac{\beta_{2}-1}{x}-\frac{2\alpha_{1}\beta_{1} x^{\beta_{1}-1}}{1-x^{2\beta_{1}}}+ \frac{2\alpha_{2}\beta_{2} x^{\beta_{2}-1}}{1-x^{2\beta_{2}}}-\frac{2\beta_{2}x^{2\beta_{2}-1}}{1-x^{2\beta_{2}}}+\frac{2\beta_{1}x^{2\beta_{1}-1}}{1-x^{2\beta_{1}}}\\
		&+\frac{\frac{2\alpha_{1}\beta_{1}p_{1} x^{\beta_{1}-1}}{1-x^{2\beta_{1}}}}{\left\{ 1-p_{1}\left(   1+\log \left(\frac{1+x^{\beta_{1}}}{1-x^{\beta_{1}}}\right)^{-\alpha_{1}}\right)\right\}} -\frac{\frac{2\alpha_{2}\beta_{2}p_{2} x^{\beta_{2}-1}}{1-x^{2\beta_{2}}}}{\left\{ 1-p_{2}\left(1+\log \left(\frac{1+x^{\beta_{2}}}{1-x^{\beta_{2}}}\right)^{-\alpha_{2}}\right)\right\}}.
	\end{align*} 
	Now, if $\beta_{1}=\beta_{2}=\beta, p_{1}=p_{2}=p$, then
	\begin{equation}\label{sto}
		= \frac{2\beta x^{\beta-1}}{1-x^{2\beta}} (\alpha_{2}-\alpha_{1})+ \frac{2\beta p x^{\beta-1}}{1-x^{2\beta}} \left\{\frac{\alpha_{1}}{ \left\{ 1-p\left(1+\log \left(\frac{1+x^{\beta}}{1-x^{\beta}}\right)^{-\alpha_{1}}\right)\right\} }- \frac{\alpha_{2}}{\left\{ 1-p\left(1+\log \left(\frac{1+x^{\beta}}{1-x^{\beta}}\right)^{-\alpha_{2}}\right)\right\}}  \right\}.
	\end{equation} 
	Now, we can see from equation \eqref{sto}, if $\beta_{1}=\beta_{2}=\beta, p_{1}=p_{2}=p$, and $\alpha_{1}>\alpha_{2}$, then $\frac{d}{dx}\log\left\{\frac{f_{X}(x)}{f_{Y}(x)}\right\}\leq 0$  implies that $X\leq_{lr}Y$. Consequently, the other relations also holds.
	\subsection{Order statistics} 
	Let us consider a random sample $X_1, X_2, \ldots, X_n$ taken from the RTUOMG distribution, with corresponding order statistics $X_{(1)}, X_{(2)}, \ldots, X_{(n)}$. The PDF of the $r^{th}$ order statistic $X_{(r)}$ (where $r=1,2,\ldots,n$) is given by 
	\begin{equation}\label{eq:rtpdf}
		f_{X_{(r)}}(x)= \frac{f(x)}{B(r,n-r+1)} \sum_{i=0}^{r-1} \binom{r-1}{i} (-1)^i \left(1-F(x)\right)^{n+i-r}.
	\end{equation}
	The CDF of $r^{th}$ order statistics ($X_{r}$) are given by (see \cite{singh2020induced})
	\begin{equation*}
		F_{X_{(r)}}(x)= \sum_{j=r}^{n} \binom{n}{j} F^{j}(x) [ 1-F(x) ]^{n-j}. 
	\end{equation*}
	This expression can also be written as
	\begin{equation}\label{eq:rthcdf}
		F_{X_{(r)}}(x) = \sum_{j=r}^{n} \sum_{l=0}^{n-j} \binom{n}{j} \binom{n-j}{l} (-1)^l F^{j+l}(x).
	\end{equation}
	The PDF and CDF of $r^{th}$ order statistics of RTUOMG distribution are respectively, given as 
	\begin{equation}
		f_{X_{(r)}}(x)= \frac{f(x)}{B(r,n-r+1)} \sum_{i=0}^{r-1} \binom{r-1}{i} (-1)^i \left[ \omega_{1}^{(-2\alpha,\beta)}(x) -p \omega_{1}^{(-2\alpha,\beta)}(x) \log\left\{\omega_{1}^{(-2\alpha,\beta)}(x)\right\} \right]^{n+i-r}
	\end{equation}
	and 
	\begin{equation}
		F_{X_{(r)}}(x)= \sum_{j=r}^{n} \sum_{l=0}^{n-j} \binom{n}{j} \binom{n-j}{l} (-1)^{l} \left[1- \omega_{1}^{(-2\alpha,\beta)}(x)+p \omega_{1}^{(-2\alpha,\beta)}(x) \log\left\{\omega_{1}^{(-2\alpha,\beta)}(x)\right\} \right]^{j+l}. 
	\end{equation}
	In addition, the PDFs of the smallest and largest order statistics are 
	\begin{equation*}
		f_{X_{(1)}}(x)= nf(x) \left[ \omega_{1}^{(-2\alpha,\beta)}(x) -p \omega_{1}^{(-2\alpha,\beta)}(x) \log\left\{\omega_{1}^{(-2\alpha,\beta)}(x)\right\} \right]^{n-1} 
	\end{equation*}
	and
	\begin{equation*}
		f_{X_{(n)}}(x)= nf(x) \sum_{i=0}^{n-1} \binom{n-1}{i} (-1)^i \left[ \omega_{1}^{(-2\alpha,\beta)}(x)-p \omega_{1}^{(-2\alpha,\beta)}(x) \log\left\{\omega_{1}^{(-2\alpha,\beta)}(x)\right\} \right]^{i},
	\end{equation*}
	respectively.
	\subsection{Record statistics}	Let $X_1, X_2, \ldots,X_{n}$ be a sequence of random variables from the RTUOMG distribution, and let $U_1, U_2, \ldots, U_n$ and  $L_1, L_2, \ldots, L_n$ be the first n upper and lower record statistics, respectively, observed from the sequence $X_1, X_2, \ldots,X_{n}$. Then, the PDF of $n^{th}$ upper ($U_{n}$) and lower ($L_{n}$) record  statistic are respectively, given by (see \cite{sakthivel2022record})
	\begin{equation}
		f_{U_{n}}(u_{n})=  \frac{\left(-\log [1-F(u_{n})]\right)^{(n-1)}}{(n-1)!} f(u_{n}), \quad u_{n}>0
	\end{equation}
	and
	\begin{equation}
		f_{L_{n}}(l_{n})= \frac{\left(-\log [F(l_{n})]\right)^{(n-1)}}{(n-1)!} f(l_{n}), \quad l_{n}>0. 
	\end{equation}
	Furthermore, by substituting the CDF of the RTUOMG distribution, we get
	\begin{equation*}
		f_{U_{n}}(u_{n})=  \frac{\left(-\log \left[\omega_{1}^{(-2\alpha,\beta)}(u_n) -p \omega_{1}^{(-2\alpha,\beta)}(u_n) \log\omega_{1}^{(-2\alpha,\beta)}(u_n)\right]\right)^{(n-1)}}{(n-1)!} f(u_{n}), \quad u_{n}>0
	\end{equation*}
	and
	\begin{equation*}
		f_{L_{n}}(l_{n})=  \frac{\left(-\log \left[1-\omega_{1}^{(-2\alpha,\beta)}(l_n) +p \omega_{1}^{(-2\alpha,\beta)}(l_n) \log\omega_{1}^{(-2\alpha,\beta)}(l_n)\right]\right)^{(n-1)}}{(n-1)!} f(l_{n}), \quad l_{n}>0.
	\end{equation*}
	Moreover, the joint PDF of first n upper records ${\bf{R}}= (R_1, R_2, \ldots R_n)$ is given by
	\begin{equation*}
		f_{R}(r)=\prod_{j=1}^{n-1} h(r_{j}) f(r_{n}).
	\end{equation*} 
	Hence, 
	\begin{equation*}
		f_{R}(r)=\left\{\prod_{j=1}^{n-1} \left(\frac{2\alpha\beta r_{j}^{\beta-1}}{1-r_{j}^{2\beta}}\right)\left[\frac{1-p-p\log\omega_{1}^{(-2\alpha,\beta)}(r_j)}{1+p\log\omega_{1}^{(-2\alpha,\beta)}(r_j)}\right] \right\} f(r_{n}), 
	\end{equation*}
	where ${\bf r}= (r_{1},r_{2},\ldots,r_{n})$ denotes the observed value of ${\bf R}= (R_1, R_2, \ldots R_n)$ with $r_{1}<r_{2}<\ldots<r_{n}$, and $f(.)$ denotes the pdf in equation \eqref{pdf1}.
	\clearpage
	\newpage
	\section{Parameters Estimation}\label{section6}
	\noindent This section examines the estimation of unknown parameters of RTUOMG$(\alpha, \beta, p)$ distribution. Several methods of point estimation such as maximum likelihood (ML), ordinary least square (OLS), weighted least square (WLS), Cram{\'e}r-von Mises (CvM), and Anderson Darling (AD) are applied to calculate the estimators for unknown parameters of the proposed distribution.  
	\subsection{Maximum likelihood (ML) estimation}
	Let $X_1, X_2, \ldots, X_n$ be a random sample of size $n$ taken from the RTUOMG$(\alpha, \beta, p)$ distribution. Then, the log-likelihood function can be written as
	\begin{equation*}
		\begin{aligned}
			\label{eq:Likelihood}
			L(\alpha,\beta,p|\textbf{x})&= n\log 2\alpha+ n\log\beta+(\beta-1)\sum_{i=1}^{n}\log x_{i}-\sum_{i=1}^{n}\log(1-x_{i}^{2\beta})-\alpha\sum_{i=1}^{n} \log \left(\frac{1+x_{i}^{\beta}}{1-x_{i}^{\beta}}\right)\\
			&+\sum_{i=1}^{n}\log\left\{1-p+p\alpha\log \left(\frac{1+x_{i}^{\beta}}{1-x_{i}^{\beta}}\right)\right\},
		\end{aligned}
	\end{equation*}
	where $\textbf{x}=(x_{1},x_{2},\ldots,x_{n})$. The  maximum likelihood (ML) estimators, $\hat{\alpha}_{ML}, \hat{\beta}_{ML},$ and $\hat{p}_{ML}$ of $\alpha, \beta$, and $p$ are calculated by simultaneously solving the following non-linear equations.
	\begin{equation}
		\begin{aligned}\label{MLE}
			\frac{\partial L(\alpha,\beta,p|\textbf{x})}{\partial \alpha}&=\frac{n}{\alpha}-\sum_{i=1}^{n} \log \left( \frac{1+x_{i}^\beta}{1-x_{i}^\beta}\right) +p\sum_{i=1}^{n}\left[ \frac{\log \left( \frac{1+x_{i}^\beta}{1-x_{i}^\beta} \right)}{1-p+p\alpha \log \left( \frac{1+x_{i}^\beta}{1-x_{i}^\beta} \right)}\right]=0,\\
			\frac{\partial L(\alpha,\beta,p|\textbf{x})}{\partial \beta}&=\frac{n}{\beta}+\sum_{i=1}^{n}\log x_{i}+2\sum_{i=1}^{n} \frac{x_{i}^{2\beta} \log x_{i}}{1-x_{i}^{2\beta}}-2\alpha \sum_{i=1}^{n} \frac{x_{i}^\beta \log x_{i}}{1-x_{i}^{2\beta}}\\
			&+ 2p\alpha\sum_{i=0}^{n}\left[\frac{\frac{x_{i}^\beta \log x_{i}}{1-x_{i}^{2\beta}}}{1-p+p\alpha \log \left( \frac{1+x_{i}^\beta}{1-x_{i}^\beta} \right)} \right] =0,\\
			\frac{\partial L(\alpha,\beta,p|\textbf{x})}{\partial p}&= \sum_{i=1}^{n}\left[\frac{\alpha\log\left(\frac{1+x_{i}^{\beta}}{1-x_{i}^{\beta}}\right)-1}{1-p+p\alpha \log \left( \frac{1+x_{i}^\beta}{1-x_{i}^\beta} \right)}\right]=0.
		\end{aligned}
	\end{equation}
	\noindent There are several numerical techniques such as Broyden-Flecther-Goldfarb-Shanno (BFGS) and Nelder-Mead (NM) that can be use to solve the above non-linear equations obtained in (\ref{MLE}). These techniques can be easily applied using optim() function in R-programming.
	\subsection{Ordinary least squares (OLS) and weighted least squares (WLS) estimation:}
	According to \cite{swain1988least}, the OLS estimators of unknown parameters can be obtained by minimizing the sum of squares differences between the vector of uniformized order statistics and the corresponding vector of expected values. 
	Let $X_{(1)}, X_{(2)}, \ldots, X_{(n)}$ denote the order statistics of a random sample $X_1, X_2,\ldots, X_n$ of size $n$ taken from  RTUOMG distribution. Then, the mean and variance of the 
	empirical cumulative distribution function (ecdf), $F(X_{(i)})$ are as follows 
	\begin{equation*}
		E(F(X_{(i)}))=\frac{i}{n+1};  i=1,2,\ldots,n
	\end{equation*}
	and
	\begin{equation*}
		V(F(X_{(i)}))=\frac{i(n-i+1)}{(n+1)^2(n+2)};  i=1,2,\ldots,n.
	\end{equation*}
	Based on the mean and variance of $F(X_{(i)})$, we can employ the following two methods of least squares estimation.\\ 
	\noindent \textbf{Ordinary least squares (OLS) estimation:} The OLS estimators, $\hat{\alpha}_{OLS}$, $\hat{\beta}_{OLS}$, and $\hat{p}_{OLS}$ of RTUOMG$(\alpha, \beta, p)$ distribution can be obtained by minimizing the following equation
	\begin{equation}
		\label{eq:LSE}
		Z(\alpha,\beta,p)=\sum_{i=1}^{n}\left( F(X_{(i)})- \frac{i}{n+1}\right)^2.
	\end{equation}

	\noindent By considering the distribution function given in equation \eqref{CDF} and equation \eqref{eq:LSE}, we can derive the following three non-linear equations by taking partial derivatives with respect to parameters  $\alpha$, $\beta$, and $p$  	
	\begin{equation}\label{le1}
		\frac{\partial}{\partial \alpha}Z(\alpha, \beta, p)=\sum_{i=1}^{n}\left(\frac{\partial F(X_{i},\alpha,\beta,p)}{\partial \alpha}\right) \left( F(X_{(i)})- \frac{i}{n+1}\right)=0,
	\end{equation}
	\begin{equation}\label{le2}
		\frac{\partial}{\partial \beta}Z(\alpha, \beta, p)=\sum_{i=1}^{n}\left(\frac{\partial F(X_{i},\alpha,\beta,p)}{\partial \beta}\right) \left( F(X_{(i)})- \frac{i}{n+1}\right)=0,
	\end{equation}
	\begin{equation}\label{lee3}
		\frac{\partial}{\partial p}Z(\alpha, \beta, p)=\sum_{i=1}^{n}\left(\frac{\partial F(X_{i},\alpha,\beta,p)}{\partial p}\right)\left( F(X_{(i)})- \frac{i}{n+1}\right)=0.
	\end{equation}
	
	\noindent These equations capture the relationship between the observed data and the parameters of interest. The solution of equations \eqref{le1}-\eqref{lee3} can be obtained by using a non-linear equation solver technique such as Broyden-Flecther-Goldfarb-Shanno (BFGS) and Nelder-Mead (NM) technique etc.\\
	\textbf{Weighted least squares (WLS) estimation:} The WLS estimator follows a similar procedure to the OLS estimator, where the objective is to minimize the weighted sum of squares differences. The WLS estimators, $\hat{\alpha}_{WLS}$, $\hat{\beta}_{WLS}$, and $\hat{p}_{WLS}$ for the unknown parameters $\alpha$, $\beta$, and $p$ can be computed by minimizing
	\begin{equation*}
		\label{eq:WLSE}
		W(\alpha,\beta,p)= \sum_{i=1}^{n} \eta_{i}\left(F(X_{(i)})- \frac{i}{n+1} \right)^2,
	\end{equation*}
	with respect to the unknown parameters.
	This function is minimized by adjusting the values of $\alpha$, $\beta$, and $p$ such that the sum of squared differences between the observed values $F(X_{(i)})$ and the corresponding expected values $\frac{i}{n+1}$ is minimized through weight  factor $\eta_{i}$.
	The factor $\eta_{i}$ is derived from the inverse of the variance of $F(X_{(i)})$ and depends on the specific distribution being considered. 
	\subsection{Cram{\'e}r-von Mises (CvM) estimation:} Cram{\'e}r-von Mises (CvM) estimator is one type of goodness-of-fit estimator and based on Cram{\'e}r-von Mises statistics (see \cite{dey2018statistical} ). This estimator is basically based on the difference between the empirical distribution function and the theoretical distribution function. 
	Therefore, the CvM estimators, $\hat{\alpha}_{CvM}$, $\hat{\beta}_{CvM}$, and $\hat{p}_{CvM}$ of unknown parameters can be obtained by minimizing the following equation with respect to parameters $\alpha, \beta$, and $p$  \\
	\begin{equation*}
		\label{eq:CVM}
		C(\alpha,\beta,p)= \frac{1}{12n}+\sum_{i=1}^{n}\left(F(X_{(i)})-\frac{2i-1}{2n} \right)^2.
	\end{equation*}
	
	\subsection{Anderson-Darling (AD) estimation:} \cite{anderson1952asymptotic} suggested an estimator based on Anderson-Darling statistic which minimizes the Anderson-Darling distance between the empirical and theoretical distribution function.
	The AD estimators, $\hat{\alpha}_{AD}$, $\hat{\beta}_{AD}$, and $\hat{p}_{AD}$  of unknown parameters can be computed by minimizing the following equation with respect to parameters $\alpha, \beta$, and $p$  \\
	\begin{equation*}
		\label{eq:ADE}
		Q(\alpha,\beta,p)= -n-\frac{1}{n}\sum_{i=1}^{n}(2i-1)\left[ \log F(X_{(i)})+\log\Bar{F}(X_{(n+1-i)})\right], 
	\end{equation*}	
	where $\Bar{F}(x)=1-F(x)$  denotes the SF of RTUOMG distribution.  For more details of this technique, one may refer to \cite{boos1982minimum} and \cite{arshad2022record}.
	\section{Simulation study}\label{section7}
	\vspace{0.2cm}
	\noindent In this section, a well-organized Monte Carlo simulation study is carried out to evaluate the performance of ML, OLS, WLS, CvM, and AD
	estimators of unknown parameters of proposed RTUOMG$(\alpha, \beta, p)$ distribution. A well known Broyden-Fletcher-Goldfarb-Shanno (BFGS) technique introduced by \cite{broyden1970convergence}, \cite{fletcher1970new}, \cite{goldfarb1970family}, and \cite{shanno1970conditioning} is used to obtain the  estimates of population parameters $\alpha, \beta,$ and $p$. This optimization technique is easily available in R (version 4.2.1) programming library.  In present simulation, we generate the data of different sample sizes $(n)$ such as $ 25, 50, 100, 150, 200$, and $250$ with different combinations of parameter values ($\alpha, \beta, p$) using inverse cumulative distribution function method. The absolute biases and mean squared errors (MSEs) of considered estimators are calculated using $500$ repetitions for each sample sizes. The functional form of these measures are as follows;
	absolute bias= $N^{-1} \sum_{i=1}^{N} |(\hat{\eta}_i-\eta)|$, mean squared errors (MSEs)= $N^{-1}\sum_{i=1}^{N} (\hat{\eta}_i-\eta)^2$,
	where $\eta=(\alpha, \beta, p)$ is the true value of the parameter, $\hat{\eta}_i=(\hat{\alpha}_i,\hat{\beta}_i,\hat{p}_i)$ is the estimated value of the parameter $\eta$ for the $i$th repeated sample, and $N$ be the number of repeated samples. 
	The main findings of this simulation study are reported in Table \ref{T2}-\ref{T5}. Table \ref{T2} presents the absolute biases and MSEs of the estimators for fix settings of parameters  $\eta=(0.3, 0.4, 0.8)$. From Table \ref{T2}, it can be observe that the absolute biases and MSEs of all the considered estimators are decreases as the sample sizes increases. Therefore, we can say that all considered estimators are consistent. The similar performance of considered estimators are obtained in case of other settings of parameter values.  One more observation, we analyze that the  ML estimation technique performs well at the estimation of small setting of parameters i.e. $\forall$ $\alpha, \beta,$ $p\in (0, 1)$ while AD estimation technique works better for the estimation of $\alpha \geq 1$ and $\beta$, $p \in (0,1)$. Moreover, in case of  $\beta \geq 1$ and $ \alpha$, $p \in (0,1)$ or $ \alpha$, $\beta \geq 1$ and $p \in (0, 0.5)$, AD estimation technique works better for the estimation of $\alpha$, $\beta$ and ML technique performs better to estimate the parameter $p$ in terms of absolute bias and MSEs. 
	The CvM technique works worse in all situations. Based on the results of the simulation study, we recommended the maximum likelihood (ML), Anderson-Darling (AD), and weighted least squares (WLS) estimation technique for the estimation of parameters of the proposed distribution. 
	\begin{table}[ht]
		\caption{ Absolute biases and MSEs of the estimators for settings of parameters $\eta=(0.3, 0.4, 0.8)$.}  
		\centering 
		\footnotesize
		\begin{tabular}{ccccccccc} 
			\hline
			\toprule 
			$N=500$&$\downarrow$Sample& \multicolumn{3}{c}{Absolute Bias} &&\multicolumn{2}{c}{\hspace{2.2cm}MSE}& \\
			\cmidrule{3-5} \cmidrule{7-9}  
			Estimators &(n)& $\hat{\alpha}$& $\hat{\beta}$ &$\hat{p}$&&$\hat{\alpha}$& $\hat{\beta}$ &$\hat{p}$  \\
			\hline
			ML&25& 0.0575402& 0.5636601& 0.3744900&& 0.0051212& 0.7029732& 0.2106772\\
			&50& 0.0460293& 0.4086689& 0.3399031 && 0.0033725& 0.2802979 & 0.1763369\\
			&100& 0.0360543& 0.3303454& 0.2885562&& 0.0021142& 0.1552839& 0.1273787\\
			&150& 0.0303397& 0.3070048& 0.2578554&& 0.0015002& 0.1224424& 0.0988758\\
			&200&0.0262740& 0.3058071& 0.2467511&& 0.0010812& 0.1132115& 0.0823941\\
			&250&0.0229363& 0.3081646& 0.2442920&& 0.0007905& 0.1094451& 0.0744624\\
			
			OLS &25&0.0684853& 0.3996601& 0.4082226&& 0.0071143& 0.4539321& 0.2519106\\
			&50&0.0579523& 0.3564157& 0.3997664&& 0.0051239& 0.2856817& 0.2362091\\
			&100&0.0425999& 0.2906173& 0.3190678&& 0.0030105& 0.1527906& 0.1524894\\
			&150&0.0334161& 0.2740858& 0.2768188&& 0.0018788& 0.1192181& 0.1120023 \\
			&200&0.0278311& 0.2632146& 0.2493752&& 0.0013178& 0.1026786& 0.0881649\\
			&250&0.0242371& 0.2575848& 0.2375982&& 0.0009229& 0.0908743& 0.0742529\\
			
			WLS &25&0.0669899& 0.4007772& 0.4041653&& 0.0067799& 0.4072548& 0.2491423\\
			&50&0.0541343& 0.3419191& 0.3724363&& 0.0045822& 0.2584859& 0.2145772\\
			&100&0.0402377& 0.2854043& 0.3030421&& 0.0027186& 0.1350776& 0.1414393\\
			&150& 0.0312444& 0.2712038& 0.2613837&& 0.0016816& 0.1098474& 0.1024181\\
			&200&0.0255624& 0.2651793& 0.2359697&& 0.0010919& 0.0973379& 0.0782198\\
			&250&0.0223360& 0.2699729& 0.2281872&& 0.0007785& 0.0926324& 0.0681823\\
			
			CvM &25&0.0636285& 0.5163060& 0.3927234&& 0.0063881&0.7009315& 0.2354711\\
			&50&0.0542338& 0.4067139& 0.3903516&& 0.0045268& 0.3618528 &0.2260682\\
			&100&0.0402142& 0.3101414& 0.3128003&& 0.0026977& 0.1739882& 0.1457470\\
			&150& 0.0315371& 0.2878147& 0.2728989&& 0.0016897& 0.1308638& 0.1080854\\
			&200&0.0264091& 0.2739417& 0.2463895&& 0.0012001& 0.1104116& 0.0858663\\
			&250&0.0230677& 0.2661871& 0.2354335&& 0.0008426& 0.0966527& 0.0726858\\
			
			AD  &25&0.0610353& 0.4128566& 0.3822879&& 0.0057994& 0.4160725& 0.2267281\\
			&50&0.0506099& 0.3458641& 0.3624216&& 0.0041280& 0.2184971& 0.2018953	\\
			&100&0.0387137& 0.2924843& 0.2995147&& 0.0024893& 0.1366311& 0.1349263\\
			&150&0.0305515& 0.2820945& 0.2624722&& 0.0015762& 0.1141870& 0.1003933 \\
			&200&0.0254312& 0.2755060& 0.2394694&& 0.0010539& 0.1018982& 0.0786142\\
			&250&0.0221263& 0.2774801& 0.2309920&& 0.0007644& 0.0963062& 0.0690289\\
			\hline 
		\end{tabular}\label{T2}
	\end{table}	
	\begin{table}[ht]
		\caption{ Absolute biases and MSEs of the estimators for settings of parameters $\eta=(1.5, 0.4, 0.8)$.}  
		\centering 
		\footnotesize
		\begin{tabular}{ccccccccc} 
			
			\hline
			$N=500$&$\downarrow$Sample& \multicolumn{3}{c}{Absolute Bias} &&\multicolumn{2}{c}{\hspace{2.2cm}MSE}& \\
			\cmidrule{3-5} \cmidrule{7-9}  
			Estimators &(n)& $\hat{\alpha}$& $\hat{\beta}$ &$\hat{p}$&&$\hat{\alpha}$& $\hat{\beta}$ &$\hat{p}$  \\
			\hline
			ML&25&0.3497333& 0.2236051& 0.5125836&& 0.19385894& 0.06425815& 0.3318941\\
			&50& 0.2597831& 0.1979366& 0.4998509&& 0.09421299& 0.04653150& 0.3204776\\
			&100&0.2012286& 0.1674090& 0.4492008&& 0.05820005& 0.03379105& 0.2773317\\
			&150& 0.1763568& 0.1426037& 0.3938933&& 0.04577178& 0.02613773& 0.2301613\\
			&200&0.1439343& 0.1229506& 0.3423611&& 0.03364347& 0.02047264& 0.1833536\\
			&250&0.1175313& 0.1062177& 0.3001589&& 0.02341726& 0.01611237& 0.1463373\\
			
			OLS&25&0.4566141& 0.1802872& 0.7527976&& 0.26447902& 0.05098879& 0.5877450\\
			&50&0.4078215& 0.1784580& 0.7482363&& 0.20304339& 0.04229995& 0.5834461\\
			&100&0.3974076& 0.1773041& 0.7411818&& 0.18009355& 0.03623652& 0.5748995\\
			&150&0.3867420& 0.1744949& 0.7320833&& 0.16853860& 0.03349793& 0.5649600 \\
			&200&0.3831334& 0.1725529& 0.7240970&& 0.16285928& 0.03189547& 0.5550598\\
			&250&0.3961742& 0.1784412& 0.7531013&& 0.16773712& 0.03318493& 0.5836874\\
			
			WLS&25&0.3860936& 0.1761375& 0.6553208&& 0.21006981 &0.04876736& 0.5005673\\
			&50&0.3308242& 0.1740738& 0.6410506&& 0.14698090& 0.03912981& 0.4807757\\
			&100&0.3087392& 0.1753061& 0.6455267&& 0.12307668& 0.03483763& 0.4748319\\
			&150& 0.3066162& 0.1774102& 0.6589561&& 0.11804590& 0.03419037& 0.4836404\\
			&200&0.3157362& 0.1813900& 0.6782137&& 0.12059840& 0.03468639& 0.4995112\\
			&250&0.3336347& 0.1888402& 0.7200899&& 0.12782737& 0.03683236& 0.5432789\\
			
			CvM&25&0.4082337& 0.2211056 &0.7438294&& 0.23055472& 0.07281502& 0.5790463\\
			&50&0.3582042& 0.1981327 &0.7386132&& 0.16498067 &0.05151793& 0.5736056\\
			&100&0.3711385& 0.1861582& 0.7329617&& 0.16097120& 0.03948357& 0.5664762\\
			&150&0.3794644& 0.1830302& 0.7413456&& 0.16140139& 0.03630204& 0.5745171\\
			&200&0.3780736& 0.1792293& 0.7321075&& 0.15716656& 0.03411290& 0.5631307\\
			&250&0.3880864& 0.1830108& 0.7549924&& 0.16130300& 0.03469556& 0.5851894\\
			
			AD&25&0.2909188& 0.1619203& 0.4006342&& 0.14709415& 0.04011143& 0.2089489\\
			&50& 0.2072535& 0.1499430& 0.4088768&& 0.06392817& 0.03124702& 0.2141149\\
			&100&0.1627509& 0.1383655& 0.3970670&& 0.03965682& 0.02430970& 0.2038295\\
			&150&0.1396941& 0.1368092& 0.4086397&& 0.02952598& 0.02346276& 0.2107771 \\
			&200&0.1340009& 0.1400162& 0.4203310&& 0.02649028& 0.02333142& 0.2161384\\
			&250&0.1246115& 0.1459661& 0.4383566&& 0.02279590& 0.02426058& 0.2238062\\
			\hline 
		\end{tabular}\label{T3}
	\end{table}	
	\begin{table}[ht]
		\caption{ Absolute biases and MSEs of the estimators for setting of parameters $\eta=(0.5,1.4, 0.7)$.} 
		\centering 
		\footnotesize
		\begin{tabular}{ccccccccc} 
			\hline
			\toprule 
			$N=500$&$\downarrow$Sample& \multicolumn{3}{c}{Absolute Bias} &&\multicolumn{2}{c}{\hspace{2.2cm}MSE}& \\
			\cmidrule{3-5} \cmidrule{7-9}  
			Estimators &(n)& $\hat{\alpha}$& $\hat{\beta}$ &$\hat{p}$&&$\hat{\alpha}$& $\hat{\beta}$ &$\hat{p}$  \\
			\hline
			ML&25&0.1120347& 0.4778450& 0.3510794&& 0.0192453& 0.4553152 &0.1653206\\
			&50&0.0882323 &0.4759241 &0.3538716 &&0.0115647 &0.3755231 &0.1699393 \\
			&100&0.0730091 &0.4901289 &0.3268296 &&0.0082510& 0.3552174& 0.1474419\\
			&150&0.0644513 &0.5374551 &0.3064261 &&0.0065860 &0.3798693 &0.1310115\\
			&200&0.0584768 &0.5666936 &0.2965667 &&0.0053447 &0.3882113 &0.1181616\\
			&250&0.0520851 &0.6055517 &0.2937361 &&0.0041196 &0.4195092 &0.1085803\\
			
			OLS&25&0.1169415 &0.5490052 &0.3767365 &&0.0207611 &0.5735533 &0.1939783\\
			&50&0.0955362 &0.4938152 &0.3775485&& 0.0136563 &0.3935783 &0.1906785\\
			&100&0.0742208 &0.4622377 &0.3252509 &&0.0087476& 0.3337550& 0.1456954\\
			&150& 0.0597227 &0.4525198 &0.2898479 &&0.0059388& 0.2942829& 0.1148964 \\
			&200& 0.0489766 &0.4411952 &0.2581814 &&0.0041646 &0.2647838 &0.0910280\\
			&250& 0.0439092 &0.4310606 &0.2446191 &&0.0032173 &0.2441803 &0.0794551\\
			
			WLS&25& 0.1137814 &0.5097649 &0.3640515&& 0.0199774 &0.4917461 &0.1856189\\
			&50&0.0897559 &0.4689411 &0.3441709&& 0.0123696 &0.3905042 &0.1703281 \\
			&100&0.0703048& 0.4475743 &0.3138069 &&0.0079280& 0.2959510& 0.1348477\\
			&150& 0.0572977& 0.4514203 &0.2805834&& 0.0056517 &0.2817592& 0.1106532 \\
			&200&0.0467773 &0.4539823 &0.2526273 &&0.0037037 &0.2667118 &0.0868828\\
			&250&0.0408495 &0.4603980 &0.2358324 &&0.0027702 &0.2676532 &0.0752283\\
			
			CvM&25&0.1153831 &0.5903472 &0.3593871 &&0.0212673 &0.7178636 &0.1758703\\
			&50&0.0933515 &0.5122568 &0.3659089 &&0.0129054 &0.4440887 &0.1778830\\
			&100&0.0710206 &0.4625675 &0.3228048 &&0.0080142& 0.3361106& 0.1407267\\
			&150&0.0573306 &0.4648073 &0.2886636 &&0.0054565 &0.3091565 &0.1126145\\
			&200& 0.0471177 &0.4544770 &0.2583327 &&0.0038578& 0.2761777 &0.0900179\\
			&250&0.0420647 &0.4416537 &0.2438570 &&0.0029645 &0.2553325 &0.0782985\\
			
			AD&25& 0.1082465 &0.4530197 &0.3458671&& 0.0182808 &0.3655762 &0.1621560\\
			&50&0.0881673 &0.4225620 &0.3487991&& 0.0118716 &0.2853728 &0.1611526 \\
			&100&0.0696419 &0.4106469 &0.3230672&& 0.0077281& 0.2279476& 0.1324356\\
			&150&0.0562887 &0.4263870 &0.3067596 &&0.0053736 &0.2224153 &0.1148221\\
			&200& 0.0489789 &0.4455409 &0.2952685 &&0.0040292 &0.2305750 &0.1036650\\
			&250&0.0434670 &0.4665140 &0.2900217&& 0.0030678 &0.2463197 &0.0952949\\
			\hline
			
		\end{tabular}\label{T4}
	\end{table}	
	\begin{table}[ht]
		\caption{ Absolute biases and MSEs of the estimators for setting of parameters $\eta=(2.5,3,0.2)$.} 
		\centering 
		\footnotesize
		\begin{tabular}{ccccccccc} 
			\hline
			\toprule 
			$N=500$&$\downarrow$Sample& \multicolumn{3}{c}{Absolute Bias} &&\multicolumn{2}{c}{\hspace{2.2cm}MSE}& \\
			\cmidrule{3-5} \cmidrule{7-9}  
			Estimators &(n)& $\hat{\alpha}$& $\hat{\beta}$ &$\hat{p}$&&$\hat{\alpha}$& $\hat{\beta}$ &$\hat{p}$  \\
			\hline
			ML&25&0.8061627& 0.4544024 &0.1484790 &&1.4487436 &0.3635223 &0.0359888\\
			&50&0.4895446 &0.3169428 &0.1301799 &&0.4510390 &0.1720310 &0.0260875\\
			&100&0.3099033& 0.2278899& 0.1101981 &&0.1644552 &0.0861656 &0.0161722\\
			&150&0.2312480& 0.1948946& 0.1026410 &&0.0930790 &0.0607583 &0.0125469\\
			&200&0.1849125& 0.1793090& 0.0922127 &&0.0583091 &0.0469182 &0.0093837\\
			&250&0.1519173& 0.1785981& 0.0936273 &&0.0380139 &0.0435949 &0.0093413\\
			
			OLS&25&0.6721036& 0.5160953& 0.1965019 &&0.9517571 &0.4045551 &0.0566855 \\
			&50& 0.4528296 &0.3629825 &0.1719837 &&0.3543042 &0.2104053 &0.0432459\\
			&100&0.2923195& 0.2512083 &0.1513627 &&0.1444090 &0.0976153 &0.0327181\\
			&150& 0.2309553& 0.1999917& 0.1387837 &&0.0850821& 0.0607940& 0.0274345\\
			&200&0.1967073 &0.1628344 &0.1302848 &&0.0587973 &0.0396424 &0.0227732 \\
			&250&0.1877359 &0.1322615 &0.1191776 &&0.0522507 &0.0278735 &0.0187148 \\
			
			WLS&25&0.6634569& 0.4632552& 0.2016422 &&0.8629090 &0.3438677 &0.0626860\\
			&50& 0.4572236 &0.3234613 &0.1733783 &&0.3559206 &0.1814369 &0.0474915 \\
			&100&0.3155410& 0.2245553 &0.1436479 &&0.1673690 &0.0811050 &0.0315608\\
			&150&  0.2360105& 0.1821378& 0.1233989&& 0.0900462 &0.0535859 &0.0214425 \\
			&200& 0.2038839& 0.1476566 &0.1091343 &&0.0651935 &0.0329917 &0.0147964\\
			&250&0.1805628& 0.1366309 &0.1015546 &&0.0500695 &0.0285151 &0.0120991\\
			
			CvM&25&0.8235136 &0.4981632& 0.1845208&& 1.6105885 &0.4277520 &0.0491581\\
			&50&0.4949517& 0.3507950 &0.1657580 &&0.4490723 &0.2131639 &0.0397025\\
			&100&0.3015803& 0.2441061 &0.1494132 &&0.1613265 &0.0950858& 0.0316030\\
			&150&0.2275366 &0.1979540 &0.1366198 &&0.0856153 &0.0597272& 0.0261653\\
			&200& 0.1919231& 0.1591142 &0.1271712&& 0.0570228 &0.0378875 &0.0212320 \\
			&250& 0.1851113 &0.1315436 &0.1184449 &&0.0511737 &0.0274468 &0.0182570\\
			
			AD&25&0.6703474& 0.4422730 &0.1777463&& 0.9582408 &0.3274166 &0.0476750 \\
			&50& 0.4453855& 0.3134355 &0.1599620 &&0.3374060 &0.1693089 &0.0382549\\
			&100&0.3046488& 0.2203649 &0.1362785 &&0.1534672 &0.0796048 &0.0256305\\
			&150&0.2329083& 0.1824005 &0.1257978 &&0.0882921 &0.0528154 &0.0203533\\
			&200&0.1950100& 0.1504189& 0.1149436 &&0.0590851 &0.0343346 &0.0154094\\
			&250&0.1768263& 0.1350550 &0.1087862 &&0.0466696 &0.0284675 &0.0129945\\
			\hline  
		\end{tabular}\label{T5}
	\end{table}	
	\clearpage
	\newpage
	\section{Data Analysis}\label{section8}
	\vspace{0.2cm}
	\noindent In this section, we have analyzed two different datasets to assess the usefulness of proposed  distribution in modeling real life data. Considered data set are also fitted to some well known distributions such as unit omega (UOMG), complementary unit gompertz (CUG), complementary unit lomax (CUL) (see \cite{guerra2021unit}), and unit weibull (UW) distribution (see \cite{mazucheli2020unit}). The probability density and distribution functions of these fitted distributions are presented in Table \ref{DENSITY}.
	\begin{table}[h]
		\caption{ Probability density and distribution function of considered distributions.} 
		\centering 
		\footnotesize
		\begin{tabular}{ccl} 
			\hline
			Distribution & Density Function & Distribution Function \\
			\hline
			UOMG& $\frac{2\alpha\beta x^{\beta-1}}{(1-x^{2\beta})}\biggr(\frac{1+x^\beta}{1-x^\beta}\biggr)^{-\alpha}$
			&$1-\left(\frac{1+x^\beta}{1-x^\beta}\right)^{-\alpha}$; $\alpha>0$, $\beta>0$, $x\in (0,1)$   \\
			&&\\
			\hline
			CUG&$\frac{\beta \log 2}{(1-\mu)^{-\beta}-1}(1-x)^{-(\beta+1)} 2^{\left(\frac{(1-x)^{-\beta}-1}{1-(1-\mu)^{-\beta}}\right)}$&$1-2^{\left(\frac{(1-x)^{-\beta}-1}{1-(1-\mu)^{-\beta}}\right)}$; $\beta>0$, $\mu \in (0,1)$, $x\in (0,1)$\\
			&&\\
			\hline
			CUL& $\frac{\log 2}{\beta(1-x)}\left[ \log \left( 1-\beta^{-1} \log(1-\mu)    \right)\right]^{-1} $&$1-\left[1-\beta^{-1}\log(1-x)\right]^{\left(\frac{-\log 2}{\log [1-\beta^{-1}\log(1-\mu)]}\right)}$;  \\
			& $\times \left[1-\beta^{-1}\log(1-x)\right]^{\left(\frac{-\log 2}{\log [1-\beta^{-1}\log(1-\mu)]}\right)-1}$& $\beta>0$, $\mu \in (0,1)$, $x\in (0,1)$\\
			&&\\
			\hline
			UW&$\frac{1}{x}\alpha\beta (-\log x)^{\beta-1}e^{-\alpha(-\log x)^{\beta}}$ &$e^{-\alpha(-\log x)^{\beta}}$; $\alpha>0$, $\beta>0$, $x\in (0,1)$\\ 
			&&\\
			\hline	
		\end{tabular}\label{DENSITY}	
	\end{table}	
	First of all, we have performed the exploratory data analysis and obtain the MLEs of the RTUOMG distribution with four other taken distributions. Then, we have done a comparative study of RTUOMG distribution with UOMG, CUG, CUL, and UW distribution. To determine the performance and appropriateness of fitted distributions, we use some selection statistics such as the value of -2log-likelihood function of fitted distribution, Akaike's information  criteria (AIC) of fitted model and some goodness of fit test statistics; Kolmogorov-Smirnov (KS), Cramer von Mises (CvM), and Anderson-Darling (AD) along with their p-values. The goodness of fit statistic with smaller value and highest p-value gives the better fit of distribution. The functional form of these measures and their performing algorithm are easily available in R-programming library.\\
	\noindent \textbf{Data Analysis 1:} Here, we consider a bladder cancer data (see \cite{lee2003statistical}) about the remission times (in months) of 128 patients. 
	To avoid the biased results, first, we have performed the outlier analysis using boxplot technique and extract one extreme observation 79.05* from data set. Since, the aim of this analysis is to verify the usefulness of proposed distribution with support $(0,1)$. Therefore, the final data set of 127 observations divided by 100 are as follows; \{0.0008, 0.0209, 0.0348, 0.0487, 0.0694, 0.0866, 0.1311, 0.2363, 0.0020, 0.0223, 0.0352, 0.0498, 0.0697, 0.0902,
	0.1329, 0.0040, 0.0226, 0.0357, 0.0506, 0.0709, 0.0922,
	0.1380, 0.2574, 0.0050, 0.0246, 0.0364, 0.0509, 0.0726,
	0.0947, 0.1424, 0.2582, 0.0051, 0.0254, 0.0370, 0.0517,
	0.0728, 0.0974, 0.1476, 0.2631, 0.0081, 0.0262, 0.0382,
	0.0532, 0.0732, 0.1006, 0.1477, 0.3215, 0.0264, 0.0388,
	0.0532, 0.0739, 0.1034, 0.1483, 0.3426, 0.0090, 0.0269,
	0.0418, 0.0534, 0.0759, 0.1066, 0.1596, 0.3666, 0.0105,
	0.0269, 0.0423, 0.0541, 0.0762, 0.1075, 0.1662, 0.4301,
	0.0119, 0.0275, 0.0426, 0.0541, 0.0763, 0.1712, 0.4612,
	0.0126, 0.0283, 0.0433, 0.0549, 0.0766, 0.1125, 0.1714,
	0.0135, 0.0287, 0.0562, 0.0787, 0.1164, 0.1736, 0.0140,
	0.0302, 0.0434, 0.0571, 0.0793, 0.1179, 0.1810, 0.0146,
	0.0440, 0.0585, 0.0826, 0.1198, 0.1913, 0.0176, 0.0325,
	0.0450, 0.0625, 0.0837, 0.1202, 0.0202, 0.0331, 0.0451,
	0.0654, 0.0853, 0.1203, 0.2028, 0.0202, 0.0336, 0.0676,
	0.1207, 0.2173, 0.0207, 0.0336, 0.0693, 0.0865, 0.1263,
	0.2269\}.
	The basic information of bladder cancer data are presented in Table \ref{cancerdata}. The reported value of skewness and kurtosis indicate that data is positively skewed with high kurtosis. For modelling this positive skewed and high kurtosis data a distribution defined on support (0,1) is needed. \\ 	
	
	\begin{table}[h]
		\caption{ The descriptive statistic for bladder cancer data set.} 
		\centering 
		\footnotesize
		\begin{tabular}{cccccccc} 
			\hline
			Minimum & 1st Quartile&  Median&    Mean& 3rd Quartile&  Maximum&Skewness&Kurtosis   \\
			0.00080& 0.03335& 0.06250& 0.08817& 0.11715& 0.46120&2.08009&5.09506\\
			\hline	
		\end{tabular}\label{cancerdata}
	\end{table}		
	\begin{figure}[h]
		\centering
		\subfloat{\includegraphics[scale=0.60, angle=0]{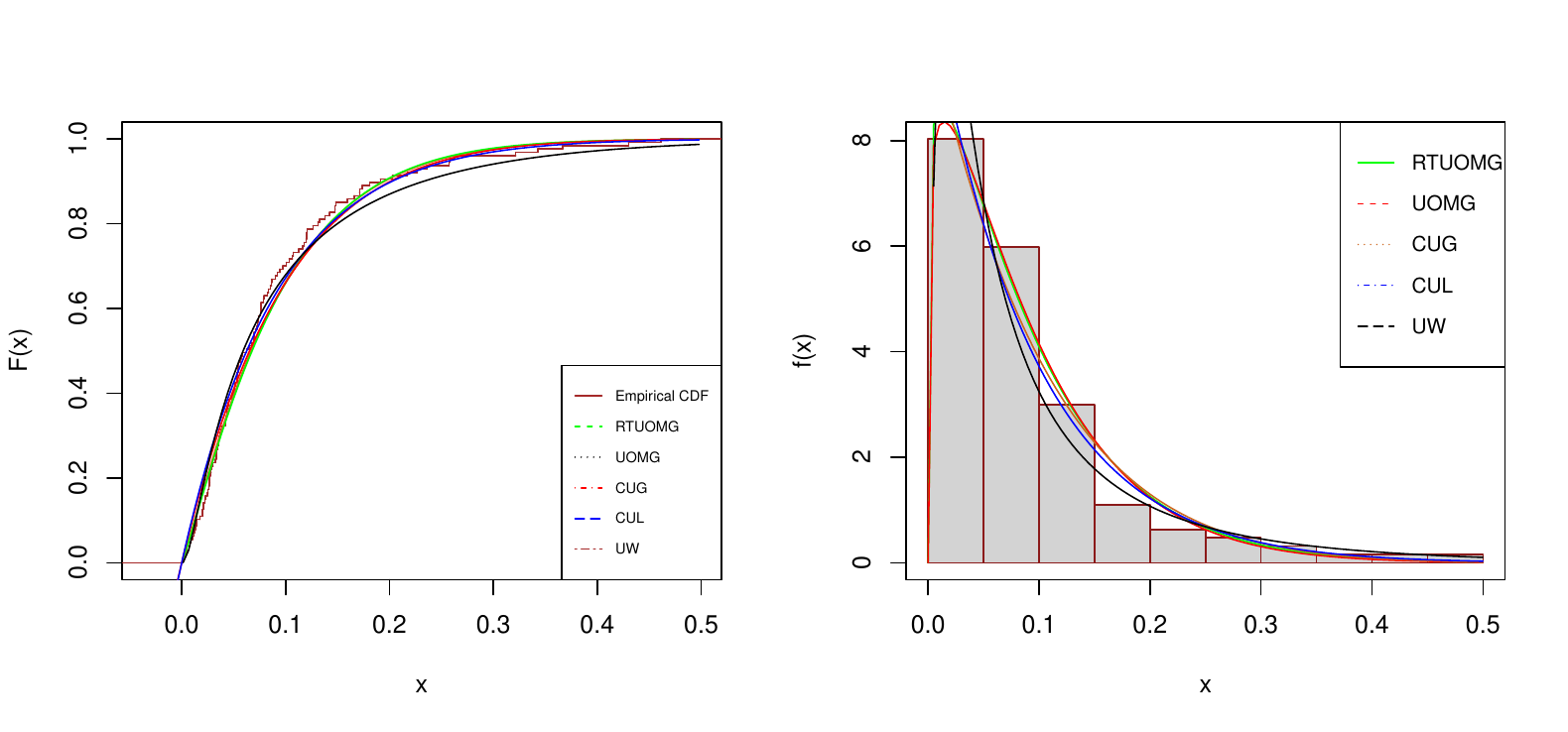}}\\
		\caption{ Fitted CDFs and PDFs for bladder cancer data.}\label{7}
	\end{figure}		
	\noindent We have fit the RTUOMG distribution with four other considered distributions for this data set and calculated MLEs with their standard errors are reported in Table \ref{MLE1}. As well as, the fitted empirical and theoretical CDFs and PDFs plots using MLEs are displayed in Figure \ref{7}. These fitted figures reveals a good fit of bladder cancer data set with RTUOMG distribution. Further, we have evaluated the values of $-2\log L$ and $AIC=2k-2\log L$, where $k$ is the number of parameters and $L$ denotes the maximized value of the likelihood function. To analyze the appropriateness of RTUOMG distribution, three other selection statistics, namely, KS, CvM, and AD along with their corresponding p-values are also calculated and reported in Table \ref{Table10}. The smallest value of AIC and goodness of fit statistics (with their high p-value) supports the best fit of RTUOMG distribution among other considered distributions.   			
	\begin{table}[h]
		\caption{ The MLEs and (standard error*) for bladder cancer data set.} 
		\centering 
		\footnotesize
		\begin{tabular}{cccc} 
			\hline
			Distribution & &MLEs (standard errors*)&  \\
			\hline
			RTUOMG& $\hat{\alpha}=$ 7.08299 (1.79777*)& $\hat{\beta}=$ 1.11156 (0.076196*)& $\hat{p}=$ 0.01 (0.19155*)\\
			UOMG&$\hat{\alpha}=$ 7.33820 (1.26419*)& $\hat{\beta}=$ 1.13426 (0.07418*)&\\
			CUG& $\hat{\mu}=$ 0.06533 (0.00660*)& $\hat{\beta}=$ 0.01000 (0.77332*)&\\
			CUL& $\hat{\mu}=$ 0.06255 (0.00646*)&  $\hat{\beta}=$ 1.42330 (2.04066*)&\\
			UW& $\hat{\alpha}=$ 0.03608 (0.00928*)& $\hat{\beta}=$ 2.84277 (0.18084*)&\\
			\hline	
		\end{tabular}\label{MLE1}
	\end{table}

	\begin{table}[h]
		\caption{ Model selection statistics for bladder cancer data set.} 
		\centering 
		\footnotesize
		\begin{tabular}{ccccccccccc} 
			
			\hline
			\toprule 
			$n=127$&& \multicolumn{3}{c}{\hspace{1.6cm}KS} &&\multicolumn{2}{c}{CvM}&& \multicolumn{2}{c}{AD}  \\
			\cmidrule{4-5} \cmidrule{7-8} \cmidrule{10-11} 
			Model &$-2\log L$ & AIC&  Statistic &P-value&& Statistic &P-value&& Statistic &P-value  \\
			\hline	
			RTUOMG&-365.0942& -359.0942& 0.0655& 0.6468&& 0.1077& 0.5492&& 0.6613& 0.5917\\
			UOMG&-362.3548& -358.3548& 0.0693& 0.5759&& 0.1145& 0.5188&& 0.6694& 0.5846\\
			CUG  &-362.3154& -358.3154& 0.0786& 0.4121&& 0.1557& 0.3730&& 1.0273& 0.3431\\
			CUL&-362.9674& -358.9674& 0.0891& 0.2654&& 0.1701& 0.3341&& 1.1275& 0.2967\\
			UW&-358.6090& -354.6090& 0.0679& 0.6011&& 0.1546& 0.3764&& 1.0428& 0.3354\\	 	
			\hline 
		\end{tabular}\label{Table10}
	\end{table}	
	
\begin{figure}[h]
		\centering
		\subfloat{\includegraphics[scale=0.60, angle=0]{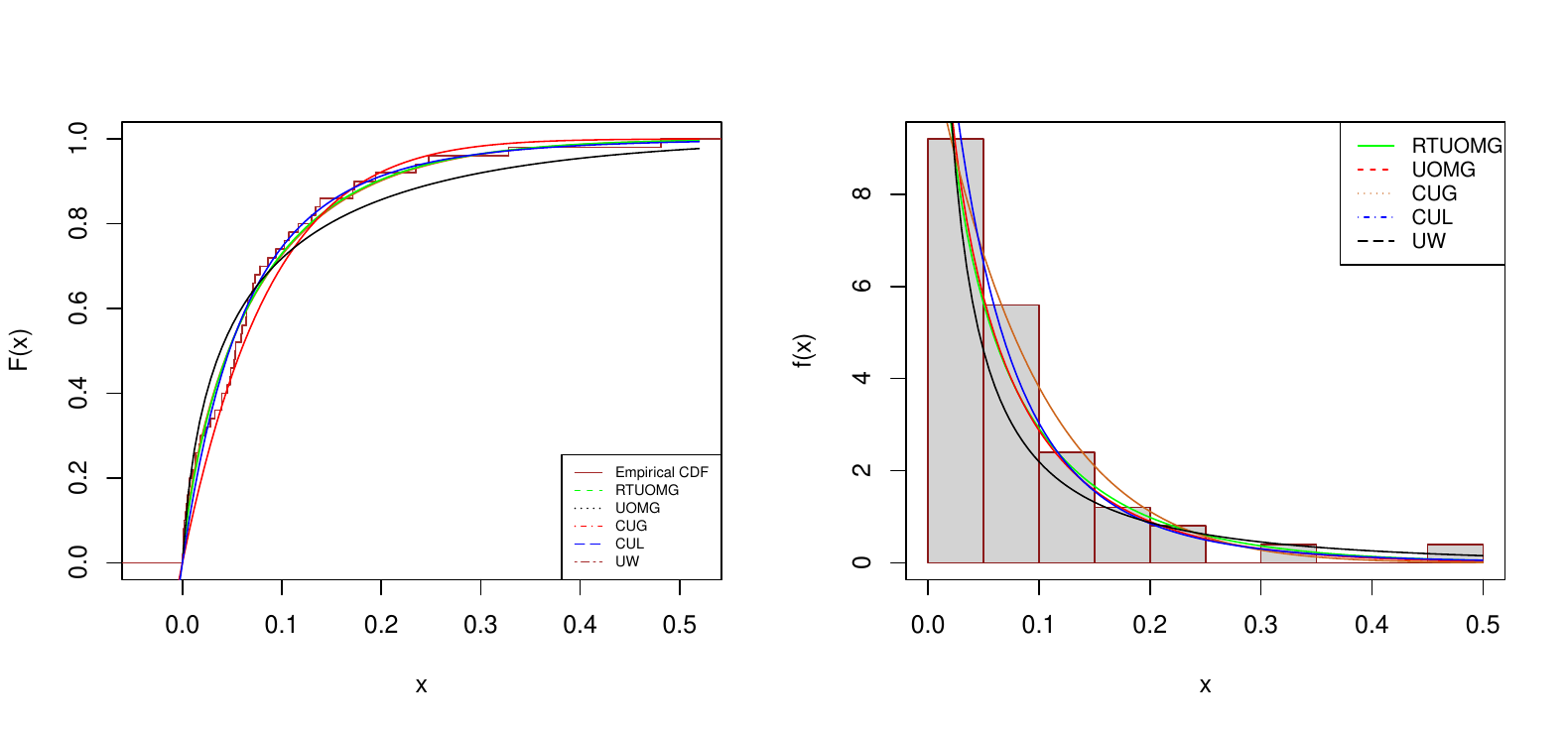}}\\
		\caption{ Fitted CDFs and PDFs for failure times data.}
		\label{Fig8}
	\end{figure}

	\noindent \textbf{Data Analysis 2:} Here, we consider a failure times (in weeks) data set of 50 components (see \cite{tanics2022record}). The values of data set divided by 100 are as follows; \{0.00013, 0.00065, 0.00111, 0.00111, 0.00163, 0.00309, 0.00426, 0.00535,
	0.00684, 0.00747, 0.00997, 0.01284, 0.01304, 0.01647, 0.01829, 0.02336, 0.02838, 0.03269 0.03977, 0.03981, 0.04520, 0.04789, 0.04849, 0.05202, 0.05291, 0.05349 0.05911, 0.06018, 0.06427, 0.06456, 0.06572, 0.07023, 0.07087, 0.07291 0.07787, 0.08596, 0.09388, 0.10261, 0.10713, 0.11658, 0.13006, 0.13388, 0.13842, 0.17152, 0.17283, 0.19418, 0.23471, 0.24777,
	0.32795, 0.48105\}. The descriptive statistics of failure times data set are given in Table \ref{Table11}. The MLEs estimators with their standard errors and model selection statistics fitted to the failure times data are reported in Table \ref{Table12}-\ref{Table13}. Figure \ref{Fig8} shows the graph of fitted CDFs and PDFs obtained from the estimates in Table \ref{Table12}. The RTUOMG distribution fits this data better than the UOMG, CUG, CUL, and UW distribution in terms of lowest values of AIC and goodness of fit statistics (with their high p-values).
	\begin{table}[h]
		\caption{ The descriptive statistic for failure times data set.} 
		\centering 
		\footnotesize
		\begin{tabular}{cccccccc} 
			\hline
			Minimum & 1st Quartile&  Median&    Mean& 3rd Quartile&  Maximum&Skewness&Kurtosis   \\
			0.00013& 0.01390& 0.05320& 0.07821& 0.10043& 0.48105 &2.37799&7.22886\\
			\hline	
		\end{tabular}\label{Table11}
	\end{table}		
	\begin{table}[h]
		\caption{ The MLEs and (standard error*) for failure times data set.} 
		\centering
		\footnotesize
		\begin{tabular}{cccc} 
			\hline
			Distribution & &MLEs (standard errors*)&  \\
			\hline
			RTUOMG& $\hat{\alpha}=$ 4.84868 (1.11256*)& $\hat{\beta}=$ 0.687588 (0.15471*)& $\hat{p}=$ 0.526608 (0.422478*)\\
			UOMG&$\hat{\alpha}=$ 4.27178 (1.03678*)& $\hat{\beta}=$ 0.818416 (0.08980*)&\\
			CUG& $\hat{\mu}=$ 0.05885 (0.00965*)& $\hat{\beta}=$ 0.01000 (1.03830*)&\\
			CUL& $\hat{\mu}=$ 0.04708 (0.00915*)&  $\hat{\beta}=$ 0.21444 (0.18292*)&\\
			UW& $\hat{\alpha}=$ 0.05624 (0.02112*)& $\hat{\beta}=$ 2.12460 (0.22026*)&\\
			\hline	
		\end{tabular}\label{Table12}
	\end{table}	
	\begin{table}[]
		\caption{ Model selection statistics for failure times data set.} 
		\centering
		\footnotesize
		\begin{tabular}{ccccccccccc} 
			
			\hline
			\toprule 
			$n=50$&& \multicolumn{3}{c}{\hspace{1.6cm}KS} &&\multicolumn{2}{c}{CvM}&& \multicolumn{2}{c}{AD}  \\
			\cmidrule{4-5} \cmidrule{7-8} \cmidrule{10-11} 
			Model &$-2\log L$ & AIC&  Statistic &P-value&& Statistic &P-value&& Statistic &P-value  \\
			\hline
			
			RTUOMG&-159.8693& -153.8693& 0.1003& 0.6960&& 0.0638& 0.7920&&  0.3342& 0.9100\\
			UOMG&-156.5532& -152.5532& 0.0973& 0.7306&& 0.0676& 0.7684&&  0.4231& 0.8248\\
			CUG&-152.4296& -148.4296& 0.1207& 0.4597&& 0.1702& 0.3343&&  1.6746& 0.1399\\
			CUL &-156.3623& -152.3623& 0.0892& 0.8212&& 0.0889& 0.6443&&  0.8247& 0.4629\\
			UW&-154.2329& -150.2329& 0.1483& 0.2214&& 0.1817& 0.3064&&  0.9602& 0.3784\\

			\hline 
		\end{tabular}\label{Table13}
	\end{table}
	\section{Conclusion}
	\noindent In this paper, we have proposed a new record-based transmuted generalization of the unit omega distribution and call it RTUOMG distribution. We derived mathematical expressions for the various important statistical quantities and graphically demonstrated the behavior of the density and hazard function. We adopted the five 
	different types of estimation techniques for estimating the unknown parameters of the proposed distribution. Using the Monte Carlo algorithm, a well organized simulation study has been performed to understand the behavior of the considered estimators for the  RTUOMG distribution with respect to absolute biases and MSEs under different setups of parameters. Based on simulation results, we recommend to use ML, AD, and WLS estimation technique for estimating the parameters of RTUOMG distribution. Moreover, the utility of the RTUOMG distribution has been also compared with some other considered distributions by analyzing two real life data sets. A dominating nature of RTUOMG distribution shows the more effectiveness and flexibility in real life modeling. \\\\
	
\noindent\textbf{Authors’ Contributions:} All authors contributed equally in this paper.\\
\noindent\textbf{Funding:} No funding.\\
\noindent\textbf{Availability of Data and Materials:} Not applicable.\\

\noindent\textbf{Declarations:}\\

\noindent\textbf{Ethical Approval:} Not required.\\
\noindent\textbf{Competing Interests:} Not applicable.\\
\noindent\textbf{Disclosure statement:}
	No potential conflict of interest was reported by the author(s).
	
	
	\newpage
	\bibliographystyle{apalike}
	\bibliography{akpomega1}
\end{document}